\documentclass[11pt,reqno]{amsart}
\usepackage{a4wide}

\numberwithin{equation}{section}
\usepackage{mathrsfs}
\usepackage{amsfonts}
\usepackage{amsmath}\allowdisplaybreaks[4]
\usepackage{stmaryrd}
\usepackage{amssymb}
\usepackage{amsthm}
\usepackage{mathrsfs}
\usepackage{url}
\usepackage{amsfonts}
\usepackage{amscd}
\usepackage{indentfirst}
\usepackage{enumerate}
\usepackage{amsmath,amsfonts,amssymb,amsthm}
\usepackage{amsmath,amssymb,amsthm,amscd}
\usepackage{graphicx,mathrsfs}
\usepackage{appendix}
\usepackage{color}

\usepackage[numbers,sort&compress]{natbib}

\newcommand{\R}{\mathbb{R}}

\renewcommand{\theequation}{\arabic{section}.\arabic{equation}}

\renewcommand{\Re}{\operatorname{Re}}
\renewcommand{\Im}{\operatorname{Im}}

\newtheorem{theo}{Theorem}[section]
\newtheorem{lem}[theo]{Lemma}

\newtheorem{prop}[theo]{Proposition}

\newtheorem{defi}[theo]{Definition}
\newtheorem{rema}[theo]{Remark}

\begin{document}

\title[]
{\textbf{Existence and local uniqueness of multi-peak solutions for the Chern-Simons-Schr\"{o}dinger system}}

\author{Qiaoqiao Hua, Chunhua Wang$^{\dagger}$ and Jing Yang}

\address[Qiaoqiao Hua]{School of Mathematics and Statistics, Central China Normal University, Wuhan 430079, China }
\email{hqq@mails.ccnu.edu.cn}

\address[Chunhua Wang]{School  of Mathematics and Statistics \& Hubei Key Laboratory of Mathematical Sciences,
  Central China Normal University,
  Wuhan, P. R. China }
\email{chunhuawang@ccnu.edu.cn}

\address[Jing Yang]{School of Science, Jiangsu University of Science and Technology, Zhenjiang 212003, China}
\email{yyangecho@163.com}

\thanks{$^\dagger$ Corresponding author: Chunhua Wang}

\date{}
\maketitle

\begin{abstract}
In the present paper, we consider the Chern-Simons-Schr\"{o}dinger system
\begin{equation}\nonumber
\left\{
\begin{aligned}
&-\varepsilon^{2}\Delta u+V(x)u+(A_{0}+A_{1}^{2}+A_{2}^{2})u=|u|^{p-2}u,\,\,\,\,x\in \R^{2},\\
&\partial_{1}A_{0}=A_{2}u^{2},\ \partial_{2}A_{0}=-A_{1}u^{2},\\
&\partial_{1}A_{2}-\partial_{2}A_{1}=-\frac{1}{2}|u|^{2},\ \partial_{1}A_{1}+\partial_{2}A_{2}=0,\\
\end{aligned}
\right.
\end{equation}
where $p>2,$ $\varepsilon>0$ is a parameter and $V:\mathbb{R}^{2}\rightarrow\mathbb{R}$ is a bounded continuous function. Under some mild assumptions on $V(x)$, we show the existence and local uniqueness of positive multi-peak solutions.
 Our methods mainly use the finite dimensional reduction method,  various local Pohozaev identities, blow-up analysis and the maximum principle.
 Because of the nonlocal terms involved by $A_{0},A_{1}$ and $A_{2},$ we have to obtain a series of new and technical estimates.
\end{abstract}

\section{Introduction}\label{sec1}
\setcounter{equation}{0}

In this paper, we are concerned with the static solution $(u,A_{0},A_{1},A_{2})$ of the Chern-Simons-Schr\"{o}dinger (CSS) system, which satisfies
\begin{equation}\label{eq1}
\left\{
\begin{aligned}
&-\varepsilon^{2}\Delta u+V(x)u+(A_{0}+A_{1}^{2}+A_{2}^{2})u=|u|^{p-2}u,\,\,\,\,x\in \R^{2},\\
&\partial_{1}A_{0}=A_{2}u^{2},\ \partial_{2}A_{0}=-A_{1}u^{2},\\
&\partial_{1}A_{2}-\partial_{2}A_{1}=-\frac{1}{2}|u|^{2},\ \partial_{1}A_{1}+\partial_{2}A_{2}=0,\\
\end{aligned}
\right.
\end{equation}
where $\varepsilon>0$ is a parameter and $V:\mathbb{R}^{2}\rightarrow\mathbb{R}$ is a bounded continuous function.

Problem \eqref{eq1} arises in the study of the following nonlinear Schr\"{o}dinger equation with the Chern-Simons gauge fields
\begin{equation}\label{intro1}
\left\{
\begin{aligned}
&-ihD_0\phi-h^2(D_1D_1+D_2D_2)\phi+V\phi=|\phi|^{p-2}\phi,\\
&\partial_0A_1-\partial_1A_0=-h\Im(\bar{\phi}D_2\phi),\\
&\partial_0A_2-\partial_2A_0=h\Im(\bar{\phi}D_1\phi),\\
&\partial_1A_2-\partial_2A_1=-\frac{1}{2}|\phi|^2,
\end{aligned}
\right.
\end{equation}
where $i$ denotes the imaginary unit, $V(x)$ is the external potential, $h$ is the Plank constant, $\partial_0=\frac{\partial}{\partial t},\ \partial_1=\frac{\partial}{\partial x_1},\ \partial_2=\frac{\partial}{\partial x_2}$ for $(t,x_1,x_2)\in\mathbb{R}^{1+2}$, $\phi:\mathbb{R}^{1+2}\to\mathbb{C}$ is a complex scalar field, $A_\mu:\mathbb{R}^{1+2}\to\mathbb{R}$ is a gauge field and $D_\mu=\partial_\mu+\frac{i}{h}A_\mu$ is a covariant derivative for $\mu=0,1,2$. It was proposed by Jackiw and Pi \cite{ref5,ref6,ref8} to describe the dynamics of a nonrelativistic solitary wave with three dimensional Chern-Simons gauge fields, which is important for explaining  electromagnetic phenomena of anyon physics. For the case $p=4$, finite time blow-up solutions were considered in \cite{ref12,ref13}, local and global well-posedness were studied in \cite{ref10,ref11} and scattering for small solutions was proved in \cite{ref14}.

The system \eqref{intro1}, so-called Chern-Simons-Schr\"{o}dinger (CSS) system, is invariant under the gauge transformation
\begin{equation}\nonumber
\phi\to e^{i\chi}\phi,\quad A_\mu\to A_\mu-h\partial_\mu\chi,\quad \mu=0,1,2,
\end{equation}
for any arbitrary $C^\infty$ function $\chi:\mathbb{R}^{1+2}\to\mathbb{R}$. If a solution $(\phi,A_0,A_1,A_2)$ to \eqref{intro1} is of the form
\begin{equation}\nonumber
\phi(t,x)=u(x)e^{i\omega t},\quad A_\mu(t,x)=A_\mu(x),\quad\mu=0,1,2,
\end{equation}
for some function $u:\mathbb{R}^2\to\mathbb{C}$, we call it a standing wave solution of frequency $\omega$. A standing wave with $\omega=0$ is referred to be static or time independent.

The results of standing wave solutions to \eqref{intro1} by variational methods have been investigated extensively in the literature.
Byeon, Huh and Seok in \cite{ref9} found standing wave solutions to \eqref{intro1} of a particular form
\begin{equation}\nonumber
\begin{aligned}
&\phi(t,x)=u(|x|)e^{i\omega t},\quad A_0(t,x)=h_1(|x|),\\
&A_1(t,x)=\frac{x_2}{|x|^2}h_2(|x|),\quad A_2(t,x)=-\frac{x_1}{|x|^2}h_2(|x|),
\end{aligned}
\end{equation}
where $V(x)$ is a constant, $\omega>0$ and $u,h_1,h_2$ are real valued functions on $[0,\infty)$ such that $h_2(0)=0$. Since the Palais-Smale condition might not hold for all $p>2$, they devised different minimization arguments for the case $p>4$, $p=4$ and $p\in(2,4)$. To be specific, static solutions were found only in the case $p=4$. In the case of $p>4$, they considered a minimization argument on the Nehari-Pohozaev manifold, while in the case of $p\in(2,4)$, they deduced solutions as minimizers on a $L^2$-sphere. They also proved the existence of a standing wave with a vortex point of order $N$ in \cite{ref15}. Based on this type of standing wave solutions, many efforts have been done to study existence \cite{ref17,ref20,ref22}, concentration \cite{ref21} and multiplicity \cite{ref23,ref25}.
The nonexistence of standing wave solutions was discussed in \cite{ref16} by applying the Derrick-Pohozaev type identities. By studying the global behavior of energy functional for $p\in(2,4)$, Pomponio and Ruiz in \cite{ref18} proved that  whether the functional is bounded from below or not depends on frequency $\omega$, which leads to the existence and non-existence of positive solutions. Furthermore, there have been some new achievements to \eqref{intro1} recently. For instance, Deng and Li \cite{ref35} obtained the ground state solution for (CSS) system.  Shen, Squassina and Yang \cite{ref34} considered the existence results for a class of gauged Schr\"{o}dinger equations with critical exponential growth and vanishing potentials.

It seems that there are very few results of standing wave solutions to \eqref{intro1} by singular perturbation arguments. To the best of our knowledge, the first result in this respect seems to be given by Pomponio and Ruiz \cite{ref19}. They considered \eqref{intro1} in a ball under homogeneous Dirichlet boundary conditions, and proved that there existed solutions for large values of the radius and those solutions were located close to the boundary. Moreover, Azzollini and Pomponio recently in \cite{ref24} established positive energy static solutions for the (CSS) systems by considering perturbed functional. Later, under the assumption that $V(x)$ is non-radial, Deng, Long and Yang in \cite{ref1} constructed a clustering solution of the singularly perturbed problem \eqref{eq1} by Lyapunov-Schmidt reduction method, which is valid for all $p>2$.

When the Chern-Simons terms vanish, \eqref{eq1} reduces to the perturbed Schr\"{o}dinger equation
\begin{equation}\label{xiu3}
-\varepsilon^2\Delta u+V(x)u=|u|^{p-2}u.
\end{equation}
Existence and local uniqueness of \eqref{xiu3} have been studied extensively in the literature. For instance, one can refer to \cite{ref4,ref27,ref28,ref29,ref30,ref31,ref32,ref33} and the references therein.

However, there seems to be no results on the existence and the local uniqueness of the multi-peak solutions concentrating at distinct points
to problem \eqref{eq1}. So we consider the existence and local uniqueness of multi-peak static solutions to problem \eqref{eq1}, which is concentrated at multiple distinct points.
Suppose that the external potential $V(x)$ satisfies:

$(V_{1})$ $V\in L^{\infty}(\mathbb{R}^{2})$ and $0<\inf_{\mathbb{R}^{2}}V(x)\leq\sup_{\mathbb{R}^{2}}V(x)<\infty;$

$(V_{2})$ There exist $k\ (k\geq2)$ distinct points $\{a^{1},\cdots,a^{k}\}\subset\mathbb{R}^{2}$ such that for every $1\leq i\leq k$, $V\in C^{\theta}(\overline{B}_{R_0}(a^i))$ for some $\theta\in(0,1)$, and
$$V(a^i)<V(x)\ \text{for} \ 0<|x-a^i|<r$$
holds for some $r$ with $0<r<R_0=1/2\min_{1\leq i,j\leq k,i\neq j}|a^i-a^j|.$

$(V_{3})$ There exist $m>1$ and $\eta>0$ such that for each $j=1,2,$ and $i=1,\cdots,k$, $V\in C^1(B_\eta(a^i))$ and
\begin{equation}\nonumber
\left\{
\begin{aligned}
V(x)&=V(a^i)+\sum_{j=1}^{2}b_{j,i}|x_j-a_j^i|^m+O(|x-a^i|^{m+1}),&\quad x\in B_\eta(a^i),\\
\frac{\partial V}{\partial x_j}&=mb_{j,i}|x_j-a^i_j|^{m-2}(x_j-a^i_j)+O(|x-a^i|^m),&\quad x\in B_\eta(a^i),
\end{aligned}
\right.
\end{equation}
where $a^i=(a^i_1,a^i_2)\in\mathbb{R}^2$ and $b_{j,i}\in\mathbb{R}$ with $b_{j,i}\neq0$.

To be precise, we give the definition of multi-peak solutions of equation \eqref{eq1} as usual.

\begin{defi}\label{def1}
Let $k\in\{1,2,\cdots\}$. We say that $u_{\varepsilon}$ is a k-peak solution of \eqref{eq1} if $u_{\varepsilon}$ satisfies

$(\romannumeral1)$ $u_{\varepsilon}$ has k local maximum points $y_{\varepsilon}^i\in\mathbb{R}^2,\ i=1,\cdots,k$, satisfying $y_{\varepsilon}^i\rightarrow a^i$ as $\varepsilon\rightarrow0$ and $a^i\neq a^j$ for $i\neq j$;

$(\romannumeral2)$ For any given small $\tau>0$, there exists $R\gg1$ such that
\begin{equation}\nonumber
|u_{\varepsilon}(x)|\leq\tau\ \text{for}\ x\in\mathbb{R}^2\setminus\cup_{i=1}^kB_{R\varepsilon}(y_{\varepsilon}^i);
\end{equation}

$(\romannumeral3)$ There exists $C>0$ such that
\begin{equation}\nonumber
||u_{\varepsilon}||_{\varepsilon}^2=\int_{\mathbb{R}^2}(\varepsilon^2|\nabla u_{\varepsilon}|^2+V(x)|u_{\varepsilon}|^2)\leq C\varepsilon^2.
\end{equation}
\end{defi}

We denote by $U^i$ the unique positive radial solution of
\begin{equation}\label{rr}
\left\{
\begin{aligned}
&-\Delta u+V(a^i)u=u^{p-1},\ u>0,\ x\in\mathbb{R}^2,\\
&u(0)=\max_{\mathbb{R}^2}u(x),\ u\in H^1({\mathbb{R}^2}).
\end{aligned}
\right.
\end{equation}
It is well-known that $U^i(x)=U^i(|x|)$ is non-degenerate and satisfies
\begin{equation}\nonumber
(U^i(r))'<0,\ \lim_{r\rightarrow\infty}\sqrt{r}e^rU^i(r)=C>0,\ \lim_{r\rightarrow\infty}\frac{(U^i(r))'}{U^i(r)}=-1.
\end{equation}
Denote
\begin{equation}\nonumber
\langle u,v \rangle_{\varepsilon}=\int_{\mathbb{R}^2}(\varepsilon^2\nabla u\nabla v+V(x)uv)\,dx \ \ \text{and}\ \  ||u||_{\varepsilon}^2=\langle u,u \rangle_{\varepsilon},
\end{equation}
and let
\begin{equation}\nonumber
H_{\varepsilon}=\{u\in H^1(\mathbb{R}^2):||u||_{\varepsilon}<\infty\}.
\end{equation}
Recall that we assume $a^i\neq a^j$ for $i\neq j$. Let $0<\delta<\min\{|a^i-a^j|/4:i\neq j\}$ and denote
\begin{equation}\nonumber
D_{\delta}=\bar{B}_{\delta}(a^1)\times\cdots\times\bar{B}_{\delta}(a^k).
\end{equation}
If $(y^1,\cdots,y^k)\in D_{\delta}$, then $|y^i-y^j|\geq|a^i-a^j|/2\geq2\delta$ for $i\neq j$.

Our first result is the following.

\begin{theo}\label{thm1}
Suppose that $p>2$ and $V(x)$ satisfies $(V_1)$ and $(V_2)$. Then for $\varepsilon>0$ sufficiently small, problem \eqref{eq1} has a k-peak solution $u_\varepsilon$ defined as in the Definition \ref{def1} concentrating around $a^i$ for $i=1,\cdots,k.$ Precisely, $u_\varepsilon$ is of the form
\begin{equation}
u_{\varepsilon}(x)=\sum_{i=1}^kU^i(\frac{x-y_\varepsilon^i}{\varepsilon})+\varphi_{\varepsilon}(x)
\end{equation}
with $y_\varepsilon^i\in\mathbb{R}^2$ and $\varphi_{\varepsilon}\in H_\varepsilon$ satisfying, for $i=1,\cdots,k$, as $\varepsilon\rightarrow0$,
\begin{equation}
y_\varepsilon^i\rightarrow a^i\quad \text{and}\quad
||\varphi_{\varepsilon}||_\varepsilon=o(\varepsilon).
\end{equation}
\end{theo}

Local uniqueness is an important topic in the study of elliptic partial differential equations.
For the local uniqueness of the k-peak solutions obtained in Theorem \ref{thm1}, we have the following theorem.
\begin{theo}\label{thm2}
Suppose that $p>2$ and $V(x)$ satisfies $(V_1),\ (V_2)$ and $(V_3)$. If $u_\varepsilon^{(i)}(x),\ i=1,2,$ are two k-peak solutions of \eqref{eq1} defined as in the Definition \ref{def1} concentrating around $a^i$ for $i=1,\cdots,k$, then for $\varepsilon>0$ sufficiently small, $u_\varepsilon^{(1)}(x)\equiv u_\varepsilon^{(2)}(x)$ must be of the form
\begin{equation}
u_{\varepsilon}(x)=\sum_{i=1}^kU^i(\frac{x-y_\varepsilon^i}{\varepsilon})+\varphi_{\varepsilon}(x)
\end{equation}
with $y_\varepsilon^i\in\mathbb{R}^2$ and $\varphi_{\varepsilon}\in H_\varepsilon$ satisfying, for $i=1,\cdots,k$, as $\varepsilon\rightarrow0$,
\begin{equation}
|y_\varepsilon^i-a^i|=o(\varepsilon)\quad \text{and}\quad ||\varphi_{\varepsilon}||_\varepsilon=O(\varepsilon^{1+m}).
\end{equation}
\end{theo}

\begin{rema}
In fact, we can consider the assumption $(V_{3})$ as a more general assumption $(\tilde{V}_{3})$ as follows:
$(\tilde{V}_{3})$ There exist $m_j>1$ and $\eta>0$ such that for each $j=1,2,$ and $i=1,\cdots,k$, $V\in C^1(B_\eta(a^i))$ and
\begin{equation}\nonumber
\left\{
\begin{aligned}
V(x)&=V(a^i)+\sum_{j=1}^{2}b_{j,i}|x_j-a_j^i|^{m_{j}}+O(|x-a^i|^{m_{j}+1}),&\quad x\in B_\eta(a^i),\\
\frac{\partial V}{\partial x_j}&=m_{j}b_{j,i}|x_j-a^i_j|^{m_{j}-2}(x_j-a^i_j)+O(|x-a^i|^{m_{j}}),&\quad x\in B_\eta(a^i),
\end{aligned}
\right.
\end{equation}
where $a^i=(a^i_1,a^i_2)\in\mathbb{R}^2$ and $b_{j,i}\in\mathbb{R}$ with $b_{j,i}\neq0,$
that is the expansion conditions of $V(x)$ and $\frac{\partial V(x)}{\partial x_{j}}$  having different orders in different direction. In this case, Theorem \ref{thm2}
also holds true. To prove it, one only need to make some minor modifications.

\end{rema}

\begin{rema}
To our best knowledge, this is the first time to consider the local uniqueness of concentrated solutions for problem \eqref{eq1}. 
\end{rema}

We will prove Theorem \ref{thm1} by the finite dimensional reduction method. Although it is standard (see \cite{CPY-2021}), we have to
overcome
 some difficulties caused by the nonlocal terms $A_0,\,A_1$ and $A_2,$
 which make computations more complicated than the usual Schr\"odinger equation. To prove Theorem \ref{thm2}, inspired by \cite{ref4,DLY-15-JMPA,GPY-17-PLMS,LPW-2020-CV}
  we mainly argue by contradiction, which involves local Pohozaev identities, blow-up analysis and the maximum principle.
  In the local Pohozaev identities, there are two more nonlocal terms and we have to estimate them very carefully.

Now, we give the main idea for the proof of main results. We will find solutions of problem \eqref{eq1} by looking for critical points of the associated functional
\begin{equation}\label{qq}
\begin{aligned}
J_{\varepsilon}(u,A_{0},A_{1},A_{2})={}&\frac{1}{2}\int_{\mathbb{R}^{2}}(\varepsilon^{2}|\nabla u|^{2}+V(x)|u|^2+(A_{0}+A_{1}^{2}+A_{2}^{2})|u|^2)\,dx\\
{}&+\frac{1}{2}\int_{\mathbb{R}^{2}}(A_0F_{12}+A_1\partial_2A_0-A_2\partial_1A_0)\,dx
-\frac{1}{p}\int_{\mathbb{R}^{2}}|u|^p\,dx,
\end{aligned}
\end{equation}
where $F_{12}=\partial_1A_2-\partial_2A_1.$
Precisely in \cite{ref22}, $A_0,A_1$ and $A_2$ can be expressed as functions of $u$ by \eqref{eq1} with $K_i(x)=-\frac{x_i}{2\pi|x|^2}$, i.e.
\begin{equation}
A_1=A_1(u)=\frac{1}{2}K_2\ast(|u|^2)=-\frac{1}{4\pi}\int_{\mathbb{R}^2}\frac{x_2-y_2}{|x-y|^2}|u(y)|^2\,dy,
\end{equation}
\begin{equation}
A_2=A_2(u)=-\frac{1}{2}K_1\ast(|u|^2)=\frac{1}{4\pi}\int_{\mathbb{R}^2}\frac{x_1-y_1}{|x-y|^2}|u(y)|^2\,dy,
\end{equation}
and
\begin{equation}
\begin{aligned}
A_0&=A_0(u)=K_1\ast(A_2|u|^2)-K_2\ast(A_1|u|^2)\\
&=-\frac{1}{8\pi^2}\sum_{i=1}^2\int_{\mathbb{R}^2}\frac{x_i-y_i}{|x-y|^2}|u(y)|^2\Big(
\int_{\mathbb{R}^2}\frac{y_i-z_i}{|y-z|^2}|u(z)|^2\,dz\Big)\,dy.
\end{aligned}
\end{equation}
Moreover, problem \eqref{qq} can be reduced to find critical points of the functional
\begin{equation}
\begin{aligned}
I_{\varepsilon}(u)={}&\frac{1}{2}\int_{\mathbb{R}^{2}}(\varepsilon^{2}|\nabla u|^{2}+V(x)|u|^2)\,dx-\frac{1}{p}\int_{\mathbb{R}^{2}}|u|^p\,dx\\
{}&+\frac{1}{2}\int_{\mathbb{R}^2}\Big(-\frac{1}{4\pi}\int_{\mathbb{R}^2}\frac{x_2-y_2}{|x-y|^2}u^2(y)\,dy\Big)^2u^2(x)\,dx\\
{}&+\frac{1}{2}\int_{\mathbb{R}^2}\Big(\frac{1}{4\pi}\int_{\mathbb{R}^2}\frac{x_1-y_1}{|x-y|^2}u^2(y)\,dy\Big)^2u^2(x)\,dx.\\
\end{aligned}
\end{equation}
For simplicity, letting $U_{\varepsilon,y^i}^i(x)=U^i(\frac{x-y^i}{\varepsilon})$ for $i=1,\cdots,k$, then $U_{\varepsilon,y^i}^i$ satisfies
\begin{equation}
-\varepsilon^2\Delta U_{\varepsilon,y^i}^i+V(a^i)U_{\varepsilon,y^i}^i=(U_{\varepsilon,y^i}^i)^{p-1},\,\,\,\,x\in \R^{2}.
\end{equation}
Letting $W_{\varepsilon,Y}=\sum_{i=1}^{k}U_{\varepsilon,y^i}^i$, by using the finite dimensional reduction method, we want to construct a k-peak solution to equation \eqref{eq1} of the form
\begin{equation}\label{form}
u_{\varepsilon}=W_{\varepsilon,Y_{\varepsilon}}+\varphi_{\varepsilon,Y_{\varepsilon}},
\end{equation}
where $Y_{\varepsilon}=(y_{\varepsilon}^1,\cdots,y_{\varepsilon}^k),\ y_{\varepsilon}^i\rightarrow a^i$ as $\varepsilon\rightarrow0$ for each $i=1,\cdots,k,$ and $\varphi_{\varepsilon,Y_{\varepsilon}}$ should be appropriately chosen such that $u_{\varepsilon}$ is indeed a solution to \eqref{eq1}.

For proving the local uniqueness in Theorem \ref{thm2}, we mainly use an indirect method.
Firstly we aim to get the improved estimates of $|y_\varepsilon^i-a^i|$ and $||\varphi_{\varepsilon}||_\varepsilon$ by local Pohozaev identities. Then using local Pohozaev identities, blow-up analysis and the maximum principle, we will show $||\xi_\varepsilon||_{L^\infty(\mathbb{R}^2)}=o(1)$ as $\varepsilon\to0$ where $\xi_\varepsilon={(u_\varepsilon^{(1)}-u_\varepsilon^{(2)})}/{||u_\varepsilon^{(1)}-u_\varepsilon^{(2)}||_{L^\infty(\mathbb{R}^2)}}$, which contradicts  with $||\xi_\varepsilon||_{L^\infty(\mathbb{R}^2)}=1$ obviously.


The structure of the paper is as follows. We prove Theorem \ref{thm1} in Section \ref{sec2}. In Section \ref{sec3-1},
we prove the local uniqueness of the k-peak solutions of \eqref{eq1}.
We put the estimate of the energy functional and some Pohozaev identities in Appendix \ref{seca} and \ref{secb} respectively.

\section{The existence of k-peak solutions}\label{sec2}

\subsection{ The finite dimensional reduction}
\noindent

In this subsection, we mainly do the finite dimensional reduction process.

The main result in this subsection is as follows.
\begin{prop}\label{prop1}
There exist $\varepsilon_0>0$ and $\delta_0>0$ sufficiently small such that for all $\varepsilon\in(0, \varepsilon_0)$ and $\delta\in(0,\delta_0)$, there exists a $C^1$ map $\varphi_\varepsilon:D_{\delta}\rightarrow E_{\varepsilon,Y};\ Y\mapsto\varphi_{\varepsilon,Y}$ satisfying
\begin{equation}
\Big\langle\frac{\partial \mathcal{J}_{\varepsilon}(Y,\varphi)}{\partial\varphi}\Big|_{\varphi=\varphi_{\varepsilon,Y}},\psi\Big\rangle_\varepsilon=0,\quad\forall\psi\in E_{\varepsilon,Y}.
\end{equation}
Moreover, there exists a constant $C>0$ independent of $\varepsilon$ such that
\begin{equation}
||\varphi_{\varepsilon,Y}||_{\varepsilon}\leq C||l_{\varepsilon}||\leq C\varepsilon\Big(\varepsilon^{\theta}+\sum_{i=1}^k|V(y^i)-V(a^i)|\Big).
\end{equation}
\end{prop}

 Define
\begin{equation}\nonumber
\mathcal{J}_{\varepsilon}(Y,\varphi)=I_{\varepsilon}(W_{\varepsilon,Y}+\varphi),\ Y\in\mathbb{R}^{2k},\ \varphi\in H_{\varepsilon}.
\end{equation}
The functional $\mathcal{J}_{\varepsilon}(Y,\varphi)$ is expanded as follows:
\begin{equation}\label{eq-J}
\mathcal{J}_{\varepsilon}(Y,\varphi)=\mathcal{J}_{\varepsilon}(Y,0)+\langle l_{\varepsilon},\varphi\rangle_\varepsilon+
\frac{1}{2}\langle L_{\varepsilon}\varphi,\varphi\rangle_\varepsilon
+R_{\varepsilon}(\varphi),
\end{equation}
where
\begin{align*}
\langle l_{\varepsilon},\varphi\rangle_\varepsilon={}&\sum_{i=1}^k\int_{\mathbb{R}^2}(V(x)-V(a^i))U_{\varepsilon,y^i}^i\varphi\,dx+
\int_{\mathbb{R}^2}\Big(\sum_{i=1}^k(U_{\varepsilon,y^i}^i)^{p-1}-(\sum_{i=1}^kU_{\varepsilon,y^i}^i)^{p-1}\Big)\varphi\,dx\\
{}&+\frac{1}{16\pi^2}\int_{\mathbb{R}^2}\Big(\int_{\mathbb{R}^2}\frac{x_2-y_2}{|x-y|^2}W_{\varepsilon,Y}^2(y)\,dy\Big)^2W_{\varepsilon,Y}(x)\varphi(x)\,dx\\
{}&+\frac{1}{8\pi^2}\int_{\mathbb{R}^2}\Big(\int_{\mathbb{R}^2}\frac{x_2-y_2}{|x-y|^2}W_{\varepsilon,Y}^2(y)\,dy
\cdot\int_{\mathbb{R}^2}\frac{x_2-z_2}{|x-z|^2}W_{\varepsilon,Y}(z)\varphi(z)\,dz\Big)W_{\varepsilon,Y}^2(x)\,dx\\
{}&+\frac{1}{16\pi^2}\int_{\mathbb{R}^2}\Big(\int_{\mathbb{R}^2}\frac{x_1-y_1}{|x-y|^2}W_{\varepsilon,Y}^2(y)\,dy\Big)^2W_{\varepsilon,Y}(x)\varphi(x)\,dx\\
{}&+\frac{1}{8\pi^2}\int_{\mathbb{R}^2}\Big(\int_{\mathbb{R}^2}\frac{x_1-y_1}{|x-y|^2}W_{\varepsilon,Y}^2(y)\,dy
\cdot\int_{\mathbb{R}^2}\frac{x_1-z_1}{|x-z|^2}W_{\varepsilon,Y}(z)\varphi(z)\,dz\Big)W_{\varepsilon,Y}^2(x)\,dx,
\end{align*}
\begin{align*}
\langle L_{\varepsilon}\varphi,\varphi\rangle_\varepsilon
={}&\int_{\mathbb{R}^2}(\varepsilon^2|\nabla\varphi|^2+
V(x)\varphi^2)\,dx-(p-1)\int_{\mathbb{R}^2}W_{\varepsilon,Y}^{p-2}\varphi^2\,dx\\
{}&+\frac{1}{16\pi^2}\int_{\mathbb{R}^2}\Big(\int_{\mathbb{R}^2}\frac{x_2-y_2}{|x-y|^2}W_{\varepsilon,Y}^2(y)\,dy\Big)^2\varphi^2(x)\,dx\\
{}&+\frac{1}{8\pi^2}\int_{\mathbb{R}^2}\Big(\int_{\mathbb{R}^2}\frac{x_2-y_2}{|x-y|^2}W_{\varepsilon,Y}^2(y)\,dy
\cdot\int_{\mathbb{R}^2}\frac{x_2-z_2}{|x-z|^2}\varphi^2(z)\,dz\Big)W_{\varepsilon,Y}^2(x)\,dx\\
{}&+\frac{1}{4\pi^2}\int_{\mathbb{R}^2}\Big(\int_{\mathbb{R}^2}\frac{x_2-y_2}{|x-y|^2}W_{\varepsilon,Y}(y)\varphi(y)\,dy\Big)^2W_{\varepsilon,Y}^2(x)\,dx\\
{}&+\frac{1}{2\pi^2}\int_{\mathbb{R}^2}\Big(\int_{\mathbb{R}^2}\frac{x_2-y_2}{|x-y|^2}W_{\varepsilon,Y}^2(y)\,dy
\cdot\int_{\mathbb{R}^2}\frac{x_2-z_2}{|x-z|^2}W_{\varepsilon,Y}(z)\varphi(z)\,dz\Big)W_{\varepsilon,Y}(x)\varphi(x)\,dx\\
{}&+\frac{1}{16\pi^2}\int_{\mathbb{R}^2}\Big(\int_{\mathbb{R}^2}\frac{x_1-y_1}{|x-y|^2}W_{\varepsilon,Y}^2(y)\,dy\Big)^2\varphi^2(x)\,dx\\
{}&+\frac{1}{8\pi^2}\int_{\mathbb{R}^2}\Big(\int_{\mathbb{R}^2}\frac{x_1-y_1}{|x-y|^2}W_{\varepsilon,Y}^2(y)\,dy
\cdot\int_{\mathbb{R}^2}\frac{x_1-z_1}{|x-z|^2}\varphi^2(z)\,dz\Big)W_{\varepsilon,Y}^2(x)\,dx\\
{}&+\frac{1}{4\pi^2}\int_{\mathbb{R}^2}\Big(\int_{\mathbb{R}^2}\frac{x_1-y_1}{|x-y|^2}W_{\varepsilon,Y}(y)\varphi(y)\,dy\Big)^2W_{\varepsilon,Y}^2(x)\,dx\\
{}&+\frac{1}{2\pi^2}\int_{\mathbb{R}^2}\Big(\int_{\mathbb{R}^2}\frac{x_1-y_1}{|x-y|^2}W_{\varepsilon,Y}^2(y)\,dy
\cdot\int_{\mathbb{R}^2}\frac{x_1-z_1}{|x-z|^2}W_{\varepsilon,Y}(z)\varphi(z)\,dz\Big)W_{\varepsilon,Y}(x)\varphi(x)\,dx\\
=:{}&\langle L_{1,\varepsilon}\varphi,\varphi\rangle_\varepsilon+\langle L_{2,\varepsilon}\varphi,\varphi\rangle_\varepsilon,
\end{align*}
where $\langle L_{1,\varepsilon}\varphi,\varphi\rangle_\varepsilon=\int_{\mathbb{R}^2}(\varepsilon^2|\nabla\varphi|^2+
V(x)\varphi^2)\,dx-(p-1)\int_{\mathbb{R}^2}W_{\varepsilon,Y}^{p-2}\varphi^2\,dx$ and $\langle L_{2,\varepsilon}\varphi,\varphi\rangle_\varepsilon=\langle L_{\varepsilon}\varphi,\varphi\rangle_\varepsilon-\langle L_{1,\varepsilon}\varphi,\varphi\rangle_\varepsilon$, and
\begin{align*}
R_{\varepsilon}(\varphi)={}&-\frac{1}{p}\int_{\mathbb{R}^2}\Big((W_{\varepsilon,Y}+\varphi)^p-W_{\varepsilon,Y}^p-
pW_{\varepsilon,Y}^{p-1}\varphi-\frac{p(p-1)}{2}W_{\varepsilon,Y}^{p-2}\varphi^2\Big)\,dx\\
{}&+\frac{1}{32\pi^2}\int_{\mathbb{R}^2}\Big(\int_{\mathbb{R}^2}\frac{x_2-y_2}{|x-y|^2}\varphi^2(y)\,dy\Big)^2W_{\varepsilon,Y}^2(x)\,dx\\
{}&+\frac{1}{16\pi^2}\int_{\mathbb{R}^2}\Big(\int_{\mathbb{R}^2}\frac{x_2-y_2}{|x-y|^2}\varphi^2(y)\,dy\Big)^2W_{\varepsilon,Y}(x)\varphi(x)\,dx\\
{}&+\frac{1}{32\pi^2}\int_{\mathbb{R}^2}\Big(\int_{\mathbb{R}^2}\frac{x_2-y_2}{|x-y|^2}\varphi^2(y)\,dy\Big)^2\varphi^2(x)\,dx\\
{}&+\frac{1}{8\pi^2}\int_{\mathbb{R}^2}\Big(\int_{\mathbb{R}^2}\frac{x_2-y_2}{|x-y|^2}W_{\varepsilon,Y}(y)\varphi(y)\,dy\Big)^2\varphi^2(x)\,dx\\
{}&+\frac{1}{4\pi^2}\int_{\mathbb{R}^2}\Big(\int_{\mathbb{R}^2}\frac{x_2-y_2}{|x-y|^2}W_{\varepsilon,Y}(y)\varphi(y)\,dy\Big)^2W_{\varepsilon,Y}(x)\varphi(x)\,dx\\
{}&+\frac{1}{16\pi^2}\int_{\mathbb{R}^2}\Big(\int_{\mathbb{R}^2}\frac{x_2-y_2}{|x-y|^2}W_{\varepsilon,Y}^2(y)\,dy
\cdot\int_{\mathbb{R}^2}\frac{x_2-z_2}{|x-z|^2}\varphi^2(z)\,dz\Big)\varphi^2(x)\,dx\\
{}&+\frac{1}{8\pi^2}\int_{\mathbb{R}^2}\Big(\int_{\mathbb{R}^2}\frac{x_2-y_2}{|x-y|^2}W_{\varepsilon,Y}^2(y)\,dy
\cdot\int_{\mathbb{R}^2}\frac{x_2-z_2}{|x-z|^2}\varphi^2(z)\,dz\Big)W_{\varepsilon,Y}(x)\varphi(x)\,dx\\
{}&+\frac{1}{8\pi^2}\int_{\mathbb{R}^2}\Big(\int_{\mathbb{R}^2}\frac{x_2-y_2}{|x-y|^2}W_{\varepsilon,Y}^2(y)\,dy
\cdot\int_{\mathbb{R}^2}\frac{x_2-z_2}{|x-z|^2}W_{\varepsilon,Y}(z)\varphi(z)\,dz\Big)\varphi^2(x)\,dx\\
{}&+\frac{1}{8\pi^2}\int_{\mathbb{R}^2}\Big(\int_{\mathbb{R}^2}\frac{x_2-y_2}{|x-y|^2}W_{\varepsilon,Y}(y)\varphi(y)\,dy
\cdot\int_{\mathbb{R}^2}\frac{x_2-z_2}{|x-z|^2}\varphi^2(z)\,dz\Big)W_{\varepsilon,Y}^2(x)\,dx\\
{}&+\frac{1}{4\pi^2}\int_{\mathbb{R}^2}\Big(\int_{\mathbb{R}^2}\frac{x_2-y_2}{|x-y|^2}W_{\varepsilon,Y}(y)\varphi(y)\,dy
\cdot\int_{\mathbb{R}^2}\frac{x_2-z_2}{|x-z|^2}\varphi^2(z)\,dz\Big)W_{\varepsilon,Y}(x)\varphi(x)\,dx\\
{}&+\frac{1}{8\pi^2}\int_{\mathbb{R}^2}\Big(\int_{\mathbb{R}^2}\frac{x_2-y_2}{|x-y|^2}W_{\varepsilon,Y}(y)\varphi(y)\,dy
\cdot\int_{\mathbb{R}^2}\frac{x_2-z_2}{|x-z|^2}\varphi^2(z)\,dz\Big)\varphi^2(x)\,dx\\
{}&+\frac{1}{32\pi^2}\int_{\mathbb{R}^2}\Big(\int_{\mathbb{R}^2}\frac{x_1-y_1}{|x-y|^2}\varphi^2(y)\,dy\Big)^2W_{\varepsilon,Y}^2(x)\,dx\\
{}&+\frac{1}{16\pi^2}\int_{\mathbb{R}^2}\Big(\int_{\mathbb{R}^2}\frac{x_1-y_1}{|x-y|^2}\varphi^2(y)\,dy\Big)^2W_{\varepsilon,Y}(x)\varphi(x)\,dx\\
{}&+\frac{1}{32\pi^2}\int_{\mathbb{R}^2}\Big(\int_{\mathbb{R}^2}\frac{x_1-y_1}{|x-y|^2}\varphi^2(y)\,dy\Big)^2\varphi^2(x)\,dx\\
{}&+\frac{1}{8\pi^2}\int_{\mathbb{R}^2}\Big(\int_{\mathbb{R}^2}\frac{x_1-y_1}{|x-y|^2}W_{\varepsilon,Y}(y)\varphi(y)\,dy\Big)^2\varphi^2(x)\,dx\\
{}&+\frac{1}{4\pi^2}\int_{\mathbb{R}^2}\Big(\int_{\mathbb{R}^2}\frac{x_1-y_1}{|x-y|^2}W_{\varepsilon,Y}(y)\varphi(y)\,dy\Big)^2W_{\varepsilon,Y}(x)\varphi(x)\,dx\\
{}&+\frac{1}{16\pi^2}\int_{\mathbb{R}^2}\Big(\int_{\mathbb{R}^2}\frac{x_1-y_1}{|x-y|^2}W_{\varepsilon,Y}^2(y)\,dy
\cdot\int_{\mathbb{R}^2}\frac{x_1-z_1}{|x-z|^2}\varphi^2(z)\,dz\Big)\varphi^2(x)\,dx\\
{}&+\frac{1}{8\pi^2}\int_{\mathbb{R}^2}\Big(\int_{\mathbb{R}^2}\frac{x_1-y_1}{|x-y|^2}W_{\varepsilon,Y}^2(y)\,dy
\cdot\int_{\mathbb{R}^2}\frac{x_1-z_1}{|x-z|^2}\varphi^2(z)\,dz\Big)W_{\varepsilon,Y}(x)\varphi(x)\,dx\\
{}&+\frac{1}{8\pi^2}\int_{\mathbb{R}^2}\Big(\int_{\mathbb{R}^2}\frac{x_1-y_1}{|x-y|^2}W_{\varepsilon,Y}^2(y)\,dy
\cdot\int_{\mathbb{R}^2}\frac{x_1-z_1}{|x-z|^2}W_{\varepsilon,Y}(z)\varphi(z)\,dz\Big)\varphi^2(x)\,dx\\
{}&+\frac{1}{8\pi^2}\int_{\mathbb{R}^2}\Big(\int_{\mathbb{R}^2}\frac{x_1-y_1}{|x-y|^2}W_{\varepsilon,Y}(y)\varphi(y)\,dy
\cdot\int_{\mathbb{R}^2}\frac{x_1-z_1}{|x-z|^2}\varphi^2(z)\,dz\Big)W_{\varepsilon,Y}^2(x)\,dx\\
{}&+\frac{1}{4\pi^2}\int_{\mathbb{R}^2}\Big(\int_{\mathbb{R}^2}\frac{x_1-y_1}{|x-y|^2}W_{\varepsilon,Y}(y)\varphi(y)\,dy
\cdot\int_{\mathbb{R}^2}\frac{x_1-z_1}{|x-z|^2}\varphi^2(z)\,dz\Big)W_{\varepsilon,Y}(x)\varphi(x)\,dx\\
{}&+\frac{1}{8\pi^2}\int_{\mathbb{R}^2}\Big(\int_{\mathbb{R}^2}\frac{x_1-y_1}{|x-y|^2}W_{\varepsilon,Y}(y)\varphi(y)\,dy
\cdot\int_{\mathbb{R}^2}\frac{x_1-z_1}{|x-z|^2}\varphi^2(z)\,dz\Big)\varphi^2(x)\,dx\\
=:{}&R_{1,\varepsilon}(\varphi)+R_{2,\varepsilon}(\varphi),
\end{align*}
where $R_{1,\varepsilon}(\varphi)=-\frac{1}{p}\int_{\mathbb{R}^2}((W_{\varepsilon,Y}+\varphi)^p-W_{\varepsilon,Y}^p-
pW_{\varepsilon,Y}^{p-1}\varphi-\frac{p(p-1)}{2}W_{\varepsilon,Y}^{p-2}\varphi^2)\,dx$ and $R_{2,\varepsilon}(\varphi)=R_{\varepsilon}(\varphi)-R_{1,\varepsilon}(\varphi)$.

Define
\begin{equation}
E_{\varepsilon,Y}:=\Big\{\varphi\in H_{\varepsilon}:\langle\frac{\partial U_{\varepsilon,y^i}^i}{\partial y_{l}^i},\varphi\rangle_{\varepsilon}=0,\ i=1,\cdots,k,\ l=1,2\Big\}.
\end{equation}
The following result shows that $L_\varepsilon$ is invertible when it is restricted on $E_{\varepsilon,Y},$
which plays an essential role in carrying out the reduction argument.

\begin{lem}\label{3.4}
There exist positive constants $\varepsilon_0,\delta_0$ and $\rho$, such that for any $0<\varepsilon<\varepsilon_0,\ 0<\delta<\delta_0$ and $Y\in D_{\delta}$, there holds
\begin{equation}
||L_{\varepsilon}\varphi||\geq\rho||\varphi||_{\varepsilon},\quad\forall\varphi\in E_{\varepsilon,Y}.
\end{equation}
\end{lem}

\noindent
\begin{proof}
We prove it by a contradiction argument. Suppose that there exist $\varepsilon_n\rightarrow0$, $\delta_n\rightarrow0$, $Y^n=(y^{1,n},\cdots,y^{k,n})\in D_{\delta_n}$ and $\varphi_n\in E_{\varepsilon_n,Y^n}$ such that
\begin{equation}\nonumber
\langle L_{\varepsilon_n}\varphi_n,h\rangle_{\varepsilon_n}=o_n(1)||\varphi_n||_{\varepsilon_n}||h||_{\varepsilon_n},\ \forall\ h \in E_{\varepsilon_n,Y^n}.
\end{equation}
Without loss of generality, we may assume that $||\varphi_n||_{\varepsilon_n}=\varepsilon_n$. By the proof of Lemma \ref{3.3}, we have
\begin{equation}\nonumber
|\langle L_{2,\varepsilon_n}\varphi_n,h\rangle_{\varepsilon_n}|\leq C\varepsilon_n^2||\varphi_n||_{\varepsilon_n}||h||_{\varepsilon_n}
\end{equation}
for any $h\in E_{\varepsilon_n,Y^n}$, which implies $|\langle L_{1,\varepsilon_n}\varphi_n,h\rangle_{\varepsilon_n}|=o_n(1)||\varphi_n||_{\varepsilon_n}||h||_{\varepsilon_n}$. Letting $h=\varphi_n$, then
\begin{equation}\nonumber
\int_{\mathbb{R}^2}(\varepsilon_n^2|\nabla \varphi_n|^2+V(x)\varphi_n^2)\,dx=(p-1)\int_{\mathbb{R}^2}W_{\varepsilon_n,Y^n}^{p-2}\varphi_n^2\,dx+o_n(1)||\varphi_n||_{\varepsilon_n}^2.
\end{equation}
Since
\begin{equation}\nonumber
\begin{aligned}
(p-1)\int_{\mathbb{R}^2}W_{\varepsilon_n,Y^n}^{p-2}\varphi_n^2\,dx
&\leq C\Big(\int_{\mathbb{R}^2}W_{\varepsilon_n,Y^n}^{p}\,dx\Big)^{\frac{p-2}{p}}\Big(\int_{\mathbb{R}^2}|\varphi_n|^p\,dx\Big)^{\frac{2}{p}}\\
&\leq C\varepsilon_n^{\frac{2(p-2)}{p}}(\varepsilon_n^{2-p}||\varphi_n||_{\varepsilon_n}^p)^{\frac{2}{p}}=C||\varphi_n||_{\varepsilon_n}^2,
\end{aligned}
\end{equation}
it follows that
\begin{equation}\nonumber
\int_{\mathbb{R}^2}(\varepsilon_n^2|\nabla \varphi_n|^2+V(x)\varphi_n^2)\,dx=O(||\varphi_n||_{\varepsilon_n}^2)=O(\varepsilon_n^2).
\end{equation}
Fix $i\in\{1,\cdots,k\}$ and let $\tilde{\varphi}_{n}^i(x)=\varphi_n(\varepsilon_nx+y^{i,n})$. Then
\begin{equation}\nonumber
\int_{\mathbb{R}^2}(|\nabla\tilde{\varphi}_{n}^i|^2+V(\varepsilon_nx+y^{i,n})|\tilde{\varphi}_{n}^i|^2)\,dx
=\varepsilon_n^{-2}\int_{\mathbb{R}^2}(\varepsilon_n^2|\nabla\varphi_n|^2+V(x)\varphi_n^2)\,dx=O(1),
\end{equation}
which shows that $\{\tilde{\varphi}_{n}^i\}$ is bounded in $H^1(\mathbb{R}^2)$. Thus there exist a subsequence, still denoted by $\{\tilde{\varphi}_{n}^i\}$, and $\varphi^i\in H^1(\mathbb{R}^2)$ such that
\begin{equation}\nonumber
\begin{aligned}
&\tilde{\varphi}_{n}^i\rightharpoonup\varphi^i\quad\text{weakly in}\ H^1(\mathbb{R}^2),\\
&\tilde{\varphi}_{n}^i\rightarrow\varphi^i\quad\text{in}\ L_{loc}^p(\mathbb{R}^2)\quad(2\leq p<\infty),\\
&\tilde{\varphi}_{n}^i\rightarrow\varphi^i\quad\text{a.e. in}\ \mathbb{R}^2.
\end{aligned}
\end{equation}
It then follows that $\varphi^i$ satisfies the equation
\begin{equation}\nonumber
-\Delta\varphi^i+V(y^i)\varphi^i=(p-1)(U^i)^{p-2}\varphi^i.
\end{equation}
The non-degeneracy of the solution $U^i$ gives that $\varphi^i=\sum_{i=1}^2c_i\frac{\partial U^i}{\partial x_{i}}$. Since $\varphi_n\in E_{\varepsilon_n,Y^n}$, that is, $\langle\varphi_n,\frac{\partial U_{\varepsilon_n,y^{i,n}}^i}{\partial y_{l}^i}\rangle_{\varepsilon_n}=0$, we have
\begin{equation}\nonumber
\begin{aligned}
0&=\langle\varphi_n,\frac{\partial U_{\varepsilon_n,y^{i,n}}^i}{\partial y_{l}^i}\rangle_{\varepsilon_n}\\
&=\int_{\mathbb{R}^2}\Big(\varepsilon_n^2\nabla\varphi_n\nabla\frac{\partial U_{\varepsilon_n,y^{i,n}}^i}{\partial y_{l}^i}+V(a^i)\varphi_n\frac{\partial U_{\varepsilon_n,y^{i,n}}^i}{\partial y_{l}^i}\Big)\,dx+\int_{\mathbb{R}^2}
(V(x)-V(a^i))\varphi_n\frac{\partial U_{\varepsilon_n,y^{i,n}}^i}{\partial y_{l}^i}\,dx\\
&=(p-1)\int_{\mathbb{R}^2}(U_{\varepsilon_n,y^{i,n}}^i)^{p-2}\varphi_n\frac{\partial U_{\varepsilon_n,y^{i,n}}^i}{\partial y_{l}^i}\,dx+\int_{\mathbb{R}^2}
(V(x)-V(a^i))\varphi_n\frac{\partial U_{\varepsilon_n,y^{i,n}}^i}{\partial y_{l}^i}\,dx,
\end{aligned}
\end{equation}
which implies that
\begin{equation}\nonumber
(p-1)\int_{\mathbb{R}^2}(U^i)^{p-2}\tilde{\varphi}_{n}^i\frac{\partial U^i}{\partial x_l}\,dx+\int_{\mathbb{R}^2}
(V(\varepsilon_nx+y^{i,n})-V(a^i))\tilde{\varphi}_{n}^i\frac{\partial U^i}{\partial x_l}\,dx=0.
\end{equation}
Letting $n\rightarrow\infty$, we have $\int_{\mathbb{R}^2}(U^i)^{p-2}\varphi^i\frac{\partial U^i}{\partial x_l}\,dx=0.$ Thus $c_1=c_2=0$ and then $\varphi^i\equiv0.$ Since
\begin{equation}\nonumber
\begin{aligned}
o_n(1)\varepsilon_n^2&=\langle L_{\varepsilon_n}\varphi_n,\varphi_n\rangle_{\varepsilon_n}=\langle L_{1,\varepsilon_n}\varphi_n,\varphi_n\rangle_{\varepsilon_n}+\langle L_{2,\varepsilon_n}\varphi_n,\varphi_n\rangle_{\varepsilon_n}\\
&=||\varphi_n||_{\varepsilon_n}^2-(p-1)\int_{\mathbb{R}^2}W_{\varepsilon_n,Y^n}^{p-2}\varphi_n^2\,dx
+o_n(1)||\varphi_n||_{\varepsilon_n}^2,
\end{aligned}
\end{equation}
we obtain
\begin{equation}\label{xiu2}
(p-1)\int_{\mathbb{R}^2}W_{\varepsilon_n,Y^n}^{p-2}\varphi_n^2\,dx=(1+o_n(1))\varepsilon_n^2.
\end{equation}
By taking $R$ sufficiently large and recalling that $U^i$ decays exponentially, we have
\begin{equation}\nonumber
\int_{\mathbb{R}^2\backslash B_{R\varepsilon_n}(y^{i,n})}(U_{\varepsilon_n,y^{i,n}}^i)^{p-2}\varphi_n^2\,dx=
\varepsilon_n^2\int_{\mathbb{R}^2\backslash B_{R}(0)}
(U^i)^{p-2}(\tilde{\varphi}_{n}^i)^2\,dx=o_R(1)\varepsilon_n^2.
\end{equation}
Since $\tilde{\varphi}_{n}^i\to0$ in $L^2(B_R(0))$, we get
\begin{equation}\nonumber
\int_{B_{R\varepsilon_n}(y^{i,n})}(U_{\varepsilon_n,y^{i,n}}^i)^{p-2}\varphi_n^2\,dx=\varepsilon_n^2\int_{B_{R}(0)}
(U^i)^{p-2}(\tilde{\varphi}_{n}^i)^2\,dx=o_n(1)\varepsilon_n^2.
\end{equation}
Therefore we have
\begin{equation}\nonumber
\int_{\mathbb{R}^2}(U_{\varepsilon_n,y^{i,n}}^i)^{p-2}\varphi_n^2\,dx=o_R(1)\varepsilon_n^2+o_n(1)\varepsilon_n^2,
\end{equation}
which implies that
\begin{equation}\nonumber
\int_{\mathbb{R}^2}W_{\varepsilon_n,Y^n}^{p-2}\varphi_n^2\,dx=o_R(1)\varepsilon_n^2+o_n(1)\varepsilon_n^2.
\end{equation}
This contradicts with \eqref{xiu2}.
\end{proof}

In order to prove Proposition \ref{prop1}, from \eqref{eq-J} we have to estimate $l_\varepsilon(\varphi)$
and $R_{\varepsilon}^{(i)}(\varphi)$ respectively.

\begin{lem}\label{3.1}
Suppose that $p>2$ and $V(x)$ satisfies $(V_1)$ and $(V_2)$. Then, there exists a constant $C>0$, independent of $\varepsilon$ and $\delta$, such that for any $Y\in D_{\delta}$ there holds
\begin{equation}
||l_\varepsilon||\leq C\varepsilon\Big(\varepsilon^{\theta}+\sum_{i=1}^k|(V(y^i)-V(a^i)|\Big).
\end{equation}
\end{lem}

\noindent
\begin{proof}
For any $\varphi\in H_\varepsilon$, letting $\tilde{\varphi}(x)=\varphi(\varepsilon x)$, we obtain
\begin{equation}\label{phi}
\begin{aligned}
\int_{\mathbb{R}^2}|\varphi|^p\,dx=\varepsilon^2\int_{\mathbb{R}^2}|\tilde{\varphi}|^p\,dx&\leq C\varepsilon^2\Big(\int_{\mathbb{R}^2}(|\nabla\tilde{\varphi}|^2+|\tilde{\varphi}|^2)\,dx\Big)^\frac{p}{2}\\
&=C\varepsilon^{2-p}
\Big(\int_{\mathbb{R}^2}(\varepsilon^2|\nabla\varphi|^2+|\varphi|^2)\,dx\Big)^\frac{p}{2}\\
&\leq C\varepsilon^{2-p}||\varphi||_{\varepsilon}^p.
\end{aligned}
\end{equation}
By \eqref{phi}, the assumption $(V_2)$ and H\"{o}lder inequality, we have
\begin{align*}
{}&\Big|\sum_{i=1}^k\int_{\mathbb{R}^2}(V(x)-V(a^i))U_{\varepsilon,y^i}^i\varphi\,dx\Big|\\
\leq{}&\sum_{i=1}^k\int_{\mathbb{R}^2}|(V(x)-V(y^i))U_{\varepsilon,y^i}^i\varphi|\,dx
+\sum_{i=1}^k\int_{\mathbb{R}^2}|(V(y^i)-V(a^i))U_{\varepsilon,y^i}^i\varphi|\,dx\\
\leq{}&\sum_{i=1}^k\Big(\int_{\mathbb{R}^2}|(V(x)-V(y^i))U_{\varepsilon,y^i}^i|^2\,dx\Big)^{\frac{1}{2}}||\varphi||_{\varepsilon}
+\sum_{i=1}^k\Big(\int_{\mathbb{R}^2}|(V(y^i)-V(a^i))U_{\varepsilon,y^i}^i|^2\,dx\Big)^{\frac{1}{2}}||\varphi||_{\varepsilon}\\
\leq{}&\sum_{i=1}^k\Big(C\varepsilon^2\int_{ B_{\frac{\delta}{\varepsilon}}(0)}(|\varepsilon x|^{\theta}U^i(x))^2\,dx+C\varepsilon^2\int_{\mathbb{R}^2\backslash B_{\frac{\delta}{\varepsilon}}(0)}(U^i(x))^2\,dx\Big)^{\frac{1}{2}}||\varphi||_{\varepsilon}\\
{}&+\sum_{i=1}^k\Big(\varepsilon^2\int_{\mathbb{R}^2}|(V(y^i)-V(a^i))U^i(x)|^2\,dx\Big)^{\frac{1}{2}}||\varphi||_{\varepsilon}\\
\leq{}& C\varepsilon^{1+\theta}||\varphi||_{\varepsilon}+C\varepsilon\sum_{i=1}^k|(V(y^i)-V(a^i)|||\varphi||_{\varepsilon}.
\end{align*}
On the other hand, we derive
\begin{align*}
{}&\Big|\int_{\mathbb{R}^2}\Big(\sum_{i=1}^k(U_{\varepsilon,y^i}^i)^{p-1}-(\sum_{i=1}^kU_{\varepsilon,y^i}^i)^{p-1}\Big)\varphi\,dx\Big|\\
={}&
\left\{
\begin{aligned}
&O\Big(\sum_{i\neq j}\int_{\mathbb{R}^2}(U_{\varepsilon,y^i}^i)^{\frac{p-1}{2}}(U_{\varepsilon,y^j}^j)^{\frac{p-1}{2}}\varphi\,dx\Big),&\text{if}&\ 2<p\leq3,\\
&O\Big(\sum_{i\neq j}\int_{\mathbb{R}^2}(U_{\varepsilon,y^i}^i)^{p-2}U_{\varepsilon,y^j}^j\varphi\,dx\Big),&\text{if}&\ p>3,\\
\end{aligned}
\right.\\
\leq{}&C\varepsilon \sum_{i\neq j}e^{-\min\{\frac{p-1}{2},1\}\frac{|y^i-y^j|}{\varepsilon}}||\varphi||_{\varepsilon}\\
\leq{}&Ce^{-\frac{\tau}{\varepsilon}}||\varphi||_{\varepsilon}
\end{align*}
for some $\tau>0$. Also, by \eqref{xiu1} we have
\begin{equation}\label{equa1}
\Big|\int_{\mathbb{R}^2}\frac{x_2-y_2}{|x-y|^2}W_{\varepsilon,Y}^2(y)\,dy\Big|
\leq C\sum_{i=1}^k\int_{\mathbb{R}^2}\frac{1}{|x-y|}(U_{\varepsilon,y^i}^i(y))^2\,dy\leq C\varepsilon,
\end{equation}
and
\begin{equation}\label{equa2}
\Big|\int_{\mathbb{R}^2}\frac{x_2-y_2}{|x-y|^2}W_{\varepsilon,Y}(y)\varphi(y)\,dy\Big|
\leq C\sum_{i=1}^k\Big(\int_{\mathbb{R}^2}\frac{1}{|x-y|^2}(U_{\varepsilon,y^i}^i(y))^2\,dy\Big)^{\frac{1}{2}}
\Big(\int_{\mathbb{R}^2}\varphi^2(x)\,dx\Big)^{\frac{1}{2}}\leq C||\varphi||_{\varepsilon}.
\end{equation}
Then
\begin{equation}\nonumber
\begin{aligned}
{}&\Big|\frac{1}{16\pi^2}\int_{\mathbb{R}^2}\Big(\int_{\mathbb{R}^2}\frac{x_2-y_2}{|x-y|^2}W_{\varepsilon,Y}^2(y)\,dy\Big)^2W_{\varepsilon,Y}(x)\varphi(x)\,dx\Big|\\
\leq{}& C\varepsilon^2\int_{\mathbb{R}^2}|W_{\varepsilon,Y}(x)\varphi(x)|\,dx
\leq C\varepsilon^2\Big(\int_{\mathbb{R}^2}W_{\varepsilon,Y}^2(x)\,dx\Big)^{\frac{1}{2}}
\Big(\int_{\mathbb{R}^2}\varphi^2(x)\,dx\Big)^{\frac{1}{2}}
\leq C\varepsilon^3||\varphi||_{\varepsilon},
\end{aligned}
\end{equation}
and
\begin{equation}\nonumber
\begin{aligned}
{}&\Big|\frac{1}{8\pi^2}\int_{\mathbb{R}^2}\Big(\int_{\mathbb{R}^2}\frac{x_2-y_2}{|x-y|^2}W_{\varepsilon,Y}^2(y)\,dy
\cdot\int_{\mathbb{R}^2}\frac{x_2-z_2}{|x-z|^2}W_{\varepsilon,Y}(z)\varphi(z)\,dz\Big)W_{\varepsilon,Y}^2(x)\,dx\Big|\\
\leq{}& C\varepsilon||\varphi||_{\varepsilon}\int_{\mathbb{R}^2}W_{\varepsilon,Y}^2(x)\,dx
\leq C\varepsilon^3||\varphi||_{\varepsilon}.
\end{aligned}
\end{equation}
Analogously, we also have
\begin{equation}\nonumber
\Big|\frac{1}{16\pi^2}\int_{\mathbb{R}^2}\Big(\int_{\mathbb{R}^2}\frac{x_1-y_1}{|x-y|^2}
W_{\varepsilon,Y}^2(y)\,dy\Big)^2W_{\varepsilon,Y}(x)\varphi(x)\,dx\Big|\leq C\varepsilon^3||\varphi||_{\varepsilon},
\end{equation}
and
\begin{equation}\nonumber
\Big|\frac{1}{8\pi^2}\int_{\mathbb{R}^2}\Big(\int_{\mathbb{R}^2}\frac{x_1-y_1}{|x-y|^2}W_{\varepsilon,Y}^2(y)\,dy
\cdot\int_{\mathbb{R}^2}\frac{x_1-z_1}{|x-z|^2}W_{\varepsilon,Y}(z)\varphi(z)\,dz\Big)W_{\varepsilon,Y}^2(x)\,dx\Big|
\leq C\varepsilon^3||\varphi||_{\varepsilon}.
\end{equation}
Therefore, combining above estimates the proof is completed.
\end{proof}

\begin{lem}\label{3.2}
There exists a constant $C>0$, independent of $\varepsilon$ and $\delta$, such that for $i\in\{0,\ 1,\ 2\}$, there hold
\begin{equation}
||R_{\varepsilon}^{(i)}(\varphi)||\leq C\varepsilon^{-\min\{1,p-2\}}||\varphi||_{\varepsilon}^{\min\{3-i,p-i\}}+C(\varepsilon||\varphi||_{\varepsilon}^{3-i}+
||\varphi||_{\varepsilon}^{4-i}+\varepsilon^{-1}||\varphi||_{\varepsilon}^{5-i}+\varepsilon^{-2}||\varphi||_{\varepsilon}^{6-i}).
\end{equation}
\end{lem}

\noindent
\begin{proof}
First we estimate $R_{1,\varepsilon}$.
If $2<p\leq3$, it follows from \eqref{phi} that
\begin{equation}\nonumber
|R_{1,\varepsilon}(\varphi)|\leq C\int_{\mathbb{R}^2}|\varphi|^p\,dx\leq C\varepsilon^{2-p}||\varphi||_{\varepsilon}^p,
\end{equation}
\begin{equation}\nonumber
\begin{aligned}
|\langle R_{1,\varepsilon}'(\varphi),\psi\rangle|&\leq C\int_{\mathbb{R}^2}|\varphi|^{p-1}|\psi|\,dx\\
&\leq C\Big(\int_{\mathbb{R}^2}|\varphi|^p\,dx\Big)^{\frac{p-1}{p}}\Big(\int_{\mathbb{R}^2}|\psi|^p\,dx\Big)^{\frac{1}{p}}\\
&\leq C\varepsilon^{2-p}||\varphi||_{\varepsilon}^{p-1}||\psi||_{\varepsilon},
\end{aligned}
\end{equation}
and
\begin{equation}\nonumber
\begin{aligned}
|\langle R_{1,\varepsilon}''(\varphi)\psi,\xi\rangle|&\leq C\int_{\mathbb{R}^2}|\varphi|^{p-2}|\psi||\xi|\,dx\\
&\leq C\Big(\int_{\mathbb{R}^2}|\varphi|^p\,dx\Big)^{\frac{p-2}{p}}\Big(\int_{\mathbb{R}^2}|\psi|^p\,dx\Big)^{\frac{1}{p}}
\Big(\int_{\mathbb{R}^2}|\xi|^p\,dx\Big)^{\frac{1}{p}}\\ &\leq C\varepsilon^{2-p}||\varphi||_{\varepsilon}^{p-2}||\psi||_{\varepsilon}||\xi||_{\varepsilon}.
\end{aligned}
\end{equation}
If $p>3$, we have
\begin{equation}\nonumber
\begin{aligned}
\int_{\mathbb{R}^2}W_{\varepsilon,Y}^{p-3}|\varphi|^3\,dx&\leq C\Big(\int_{\mathbb{R}^2}W_{\varepsilon,Y}^{p}\,dx\Big)^{\frac{p-3}{p}}\Big(\int_{\mathbb{R}^2}|\varphi|^{p}\,dx\Big)^{\frac{1}{p}}\\
&\leq C\sum_{i=1}^k\Big(\int_{\mathbb{R}^2}(U_{\varepsilon,y^i}^i)^{p}\,dx\Big)^{\frac{p-3}{p}}
\Big(\int_{\mathbb{R}^2}|\varphi|^{p}\,dx\Big)^{\frac{1}{p}}\\
&\leq C(\varepsilon^2)^{\frac{p-3}{p}}
\Big(\int_{\mathbb{R}^2}|\varphi|^{p}\,dx\Big)^{\frac{1}{p}}\\
&\leq C\varepsilon^{-1}||\varphi||_{\varepsilon}^3,
\end{aligned}
\end{equation}
\begin{equation}\nonumber
\begin{aligned}
\int_{\mathbb{R}^2}W_{\varepsilon,Y}^{p-3}|\varphi|^2|\psi|\,dx&\leq C\Big(\int_{\mathbb{R}^2}W_{\varepsilon,Y}^{p}\,dx\Big)^{\frac{p-3}{p}}\Big(\int_{\mathbb{R}^2}|\varphi|^{p}\,dx\Big)^{\frac{2}{p}}
\Big(\int_{\mathbb{R}^2}|\psi|^{p}\,dx\Big)^{\frac{1}{p}}\\
&\leq C(\varepsilon^2)^{\frac{p-3}{p}}
\Big(\int_{\mathbb{R}^2}|\varphi|^{p}\,dx\Big)^{\frac{2}{p}}
\Big(\int_{\mathbb{R}^2}|\psi|^{p}\,dx\Big)^{\frac{1}{p}}\\
&\leq C\varepsilon^{-1}||\varphi||_{\varepsilon}^2||\psi||_{\varepsilon},
\end{aligned}
\end{equation}
and
\begin{equation}\nonumber
\begin{aligned}
\int_{\mathbb{R}^2}W_{\varepsilon,Y}^{p-3}|\varphi||\psi||\xi|\,dx&\leq C\Big(\int_{\mathbb{R}^2}W_{\varepsilon,Y}^{p}\,dx\Big)^{\frac{p-3}{p}}\Big(\int_{\mathbb{R}^2}|\varphi|^{p}\,dx\Big)^{\frac{1}{p}}
\Big(\int_{\mathbb{R}^2}|\psi|^{p}\,dx\Big)^{\frac{1}{p}}\Big(\int_{\mathbb{R}^2}|\xi|^{p}\,dx\Big)^{\frac{1}{p}}\\
&\leq C(\varepsilon^2)^{\frac{p-3}{p}}
\Big(\int_{\mathbb{R}^2}|\varphi|^{p}\,dx\Big)^{\frac{1}{p}}
\Big(\int_{\mathbb{R}^2}|\psi|^{p}\,dx\Big)^{\frac{1}{p}}\Big(\int_{\mathbb{R}^2}|\xi|^{p}\,dx\Big)^{\frac{1}{p}}\\
&\leq C\varepsilon^{-1}||\varphi||_{\varepsilon}||\psi||_{\varepsilon}||\xi||_{\varepsilon}.
\end{aligned}
\end{equation}
These inequalities imply that
\begin{equation}\nonumber
|R_{1,\varepsilon}(\varphi)|\leq C\int_{\mathbb{R}^2}W_{\varepsilon,Y}^{p-3}|\varphi|^3\,dx+C\int_{\mathbb{R}^2}|\varphi|^p\,dx\leq C\varepsilon^{-1}||\varphi||_{\varepsilon}^3+C\varepsilon^{2-p}||\varphi||_{\varepsilon}^p,
\end{equation}
\begin{equation}\nonumber
\begin{aligned}
|\langle R_{1,\varepsilon}'(\varphi),\psi\rangle|&\leq C\int_{\mathbb{R}^2}W_{\varepsilon,Y}^{p-3}|\varphi|^2|\psi|\,dx+C\int_{\mathbb{R}^2}|\varphi|^{p-1}|\psi|\,dx\\
&\leq C\varepsilon^{-1}||\varphi||_{\varepsilon}^2||\psi||_{\varepsilon}+C\varepsilon^{2-p}
||\varphi||_{\varepsilon}^{p-1}||\psi||_{\varepsilon},
\end{aligned}
\end{equation}
and
\begin{equation}\nonumber
\begin{aligned}
|\langle R_{1,\varepsilon}''(\varphi)\psi,\xi\rangle|&\leq C\int_{\mathbb{R}^2}W_{\varepsilon,Y}^{p-3}|\varphi||\psi||\xi|\,dx
+C\int_{\mathbb{R}^2}|\varphi|^{p-2}|\psi||\xi|\,dx\\
&\leq C\varepsilon^{-1}||\varphi||_{\varepsilon}||\psi||_{\varepsilon}||\xi||_{\varepsilon}
+C\varepsilon^{2-p}||\varphi||_{\varepsilon}^{p-2}||\psi||_{\varepsilon}||\xi||_{\varepsilon}.
\end{aligned}
\end{equation}

Next, we estimate $R_{2,\varepsilon}$. A routine computation gives rise to
\begin{equation}\label{equa3}
\begin{aligned}
\Big|\int_{\mathbb{R}^2}\frac{x_2-y_2}{|x-y|^2}\varphi^2(y)\,dy\Big|&\leq \int_{B_{\varepsilon}(x)}\frac{1}{|x-y|}\varphi^2(y)\,dy+\int_{\mathbb{R}^2\setminus B_{\varepsilon}(x)}\frac{1}{|x-y|}\varphi^2(y)\,dy\\
&\leq C\Big(\int_{B_{\varepsilon}(x)}\frac{1}{|x-y|^{\frac{3}{2}}}\,dy\Big)^{\frac{2}{3}}
\Big(\int_{\mathbb{R}^2}\varphi^6(y)\,dy\Big)^{\frac{1}{3}}+\varepsilon^{-1}\int_{\mathbb{R}^2}\varphi^2(y)\,dy\\
&\leq C\varepsilon^{\frac{1}{3}}(\varepsilon^{-4}||\varphi||_{\varepsilon}^6)^{\frac{1}{3}}+
C\varepsilon^{-1}||\varphi||_{\varepsilon}^2\\
&\leq C\varepsilon^{-1}||\varphi||_{\varepsilon}^2.
\end{aligned}
\end{equation}
As the proof of Lemma \ref{3.1}, using \eqref{equa1}, \eqref{equa2} and \eqref{equa3}, we may verify
\begin{equation}\nonumber
\Big|\frac{1}{32\pi^2}\int_{\mathbb{R}^2}\Big(\int_{\mathbb{R}^2}\frac{x_2-y_2}{|x-y|^2}\varphi^2(y)\,dy\Big)^2W_{\varepsilon,Y}^2(x)\,dx\Big|\\
\leq C||\varphi||_{\varepsilon}^4,
\end{equation}
\begin{equation}\nonumber
\Big|\frac{1}{16\pi^2}\int_{\mathbb{R}^2}\Big(\int_{\mathbb{R}^2}\frac{x_2-y_2}{|x-y|^2}\varphi^2(y)\,dy\Big)^2W_{\varepsilon,Y}(x)\varphi(x)\,dx\Big|\\
\leq C\varepsilon^{-1}||\varphi||_{\varepsilon}^5,
\end{equation}
\begin{equation}\nonumber
\Big|\frac{1}{32\pi^2}\int_{\mathbb{R}^2}\Big
(\int_{\mathbb{R}^2}\frac{x_2-y_2}{|x-y|^2}\varphi^2(y)\,dy\Big)^2\varphi^2(x)\,dx\Big|\\
\leq C\varepsilon^{-2}||\varphi||_{\varepsilon}^6,
\end{equation}
\begin{equation}\nonumber
\Big|\frac{1}{8\pi^2}\int_{\mathbb{R}^2}\Big(\int_{\mathbb{R}^2}\frac{x_2-y_2}{|x-y|^2}
W_{\varepsilon,Y}(y)\varphi(y)\,dy\Big)^2\varphi^2(x)\,dx\Big|\leq C||\varphi||_{\varepsilon}^4,
\end{equation}
\begin{equation}\nonumber
\Big|\frac{1}{4\pi^2}\int_{\mathbb{R}^2}\Big(\int_{\mathbb{R}^2}\frac{x_2-y_2}{|x-y|^2}
W_{\varepsilon,Y}(y)\varphi(y)\,dy\Big)^2W_{\varepsilon,Y}(x)\varphi(x)\,dx\Big|
\leq C\varepsilon||\varphi||_{\varepsilon}^3,
\end{equation}
\begin{equation}\nonumber
\Big|\frac{1}{16\pi^2}\int_{\mathbb{R}^2}\Big(\int_{\mathbb{R}^2}\frac{x_2-y_2}{|x-y|^2}W_{\varepsilon,Y}^2(y)\,dy
\cdot\int_{\mathbb{R}^2}\frac{x_2-z_2}{|x-z|^2}\varphi^2(z)\,dz\Big)\varphi^2(x)\,dx\Big|
\leq C||\varphi||_{\varepsilon}^4,
\end{equation}
\begin{equation}\nonumber
\Big|\frac{1}{8\pi^2}\int_{\mathbb{R}^2}\Big(\int_{\mathbb{R}^2}\frac{x_2-y_2}{|x-y|^2}W_{\varepsilon,Y}^2(y)\,dy
\cdot\int_{\mathbb{R}^2}\frac{x_2-z_2}{|x-z|^2}\varphi^2(z)\,dz\Big)W_{\varepsilon,Y}(x)\varphi(x)\,dx\Big|
\leq C\varepsilon||\varphi||_{\varepsilon}^3,
\end{equation}
\begin{equation}\nonumber
\Big|\frac{1}{8\pi^2}\int_{\mathbb{R}^2}\Big(\int_{\mathbb{R}^2}\frac{x_2-y_2}{|x-y|^2}W_{\varepsilon,Y}^2(y)\,dy
\cdot\int_{\mathbb{R}^2}\frac{x_2-z_2}{|x-z|^2}W_{\varepsilon,Y}(z)\varphi(z)\,dz\Big)\varphi^2(x)\,dx\Big|
\leq C\varepsilon||\varphi||_{\varepsilon}^3,
\end{equation}
\begin{equation}\nonumber
\Big|\frac{1}{8\pi^2}\int_{\mathbb{R}^2}\Big(\int_{\mathbb{R}^2}\frac{x_2-y_2}{|x-y|^2}W_{\varepsilon,Y}(y)\varphi(y)\,dy
\cdot\int_{\mathbb{R}^2}\frac{x_2-z_2}{|x-z|^2}\varphi^2(z)\,dz\Big)W_{\varepsilon,Y}^2(x)\,dx\Big|
\leq C\varepsilon||\varphi||_{\varepsilon}^3,
\end{equation}
\begin{equation}\nonumber
\Big|\frac{1}{4\pi^2}\int_{\mathbb{R}^2}\Big(\int_{\mathbb{R}^2}\frac{x_2-y_2}{|x-y|^2}W_{\varepsilon,Y}(y)\varphi(y)\,dy
\cdot\int_{\mathbb{R}^2}\frac{x_2-z_2}{|x-z|^2}\varphi^2(z)\,dz\Big)W_{\varepsilon,Y}(x)\varphi(x)\,dx\Big|
\leq C||\varphi||_{\varepsilon}^4,
\end{equation}
and
\begin{equation}\nonumber
\Big|\frac{1}{8\pi^2}\int_{\mathbb{R}^2}\Big(\int_{\mathbb{R}^2}\frac{x_2-y_2}{|x-y|^2}W_{\varepsilon,Y}(y)\varphi(y)\,dy
\cdot\int_{\mathbb{R}^2}\frac{x_2-z_2}{|x-z|^2}\varphi^2(z)\,dz\Big)\varphi^2(x)\,dx\Big|
\leq C\varepsilon^{-1}||\varphi||_{\varepsilon}^5.
\end{equation}
The rest terms in $R_{2,\varepsilon}$ can be estimated similarly. The assertion follows by putting all these estimates together.
\end{proof}

\begin{lem}\label{3.3}
There exists a constant $C>0$, independent of $\varepsilon$, such that
\begin{equation}
||L_{\varepsilon}\varphi||\leq C||\varphi||_{\varepsilon},\quad\forall\varphi\in H_{\varepsilon}.
\end{equation}
\end{lem}

\noindent
\begin{proof}
It is obvious that $\langle L_{\varepsilon}u,v\rangle_\varepsilon$ is bi-linear for any $u,v\in H_{\varepsilon}$. Then it is sufficient to prove $L_{\varepsilon}$ is bounded. For any $u,v\in H_{\varepsilon}$, we have
\begin{equation}\nonumber
\begin{aligned}
|\langle L_{1,\varepsilon}u,v\rangle_\varepsilon|
&\leq|\langle u,v\rangle_{\varepsilon}|+(p-1)\int_{\mathbb{R}^2}W_{\varepsilon,Y}^{p-2}|u||v|\,dx\\
&\leq||u||_{\varepsilon}||v||_{\varepsilon}+(p-1)\Big(\int_{\mathbb{R}^2}W_{\varepsilon,Y}^{p}\,dx\Big)^{\frac{p-2}{p}}
\Big(\int_{\mathbb{R}^2}|u|^p\,dx\Big)^{\frac{1}{p}}\Big(\int_{\mathbb{R}^2}|v|^p\,dx\Big)^{\frac{1}{p}}\\
&\leq||u||_{\varepsilon}||v||_{\varepsilon}+C\varepsilon^{\frac{2(p-2)}{p}}(\varepsilon^{2-p}||u||_{\varepsilon}^p)^{\frac{1}{p}}
(\varepsilon^{2-p}||v||_{\varepsilon}^p)^{\frac{1}{p}}\\
&\leq C||u||_{\varepsilon}||v||_{\varepsilon}.
\end{aligned}
\end{equation}
To estimate $|\langle L_{2,\varepsilon}u,v\rangle_\varepsilon|$, there hold
\begin{equation}\nonumber
\Big|\frac{1}{16\pi^2}\int_{\mathbb{R}^2}\Big
(\int_{\mathbb{R}^2}\frac{x_2-y_2}{|x-y|^2}W_{\varepsilon,Y}^2(y)\,dy\Big)^2u(x)v(x)\,dx\Big|
\leq C\varepsilon^2\int_{\mathbb{R}^2}|uv|\,dx
\leq C\varepsilon^2||u||_{\varepsilon}||v||_{\varepsilon},
\end{equation}
\begin{equation}\nonumber
\begin{aligned}
{}&\Big|\frac{1}{8\pi^2}\int_{\mathbb{R}^2}\Big(\int_{\mathbb{R}^2}\frac{x_2-y_2}{|x-y|^2}W_{\varepsilon,Y}^2(y)\,dy
\cdot\int_{\mathbb{R}^2}\frac{x_2-z_2}{|x-z|^2}u(z)v(z)\,dz\Big)W_{\varepsilon,Y}^2(x)\,dx\Big|\\
\leq {}&C\varepsilon\int_{\mathbb{R}^2}\Big(\int_{\mathbb{R}^2}\frac{1}{|x-z|}|u(z)v(z)|\,dz\Big)W_{\varepsilon,Y}^2(x)\,dx
=C\varepsilon\int_{\mathbb{R}^2}\Big(\int_{\mathbb{R}^2}\frac{1}{|x-z|}W_{\varepsilon,Y}^2(x)\,dx\Big)|u(z)v(z)|\,dz\\
\leq{}&C\varepsilon^2\int_{\mathbb{R}^2}|uv|\,dx\leq C\varepsilon^2||u||_{\varepsilon}||v||_{\varepsilon},
\end{aligned}
\end{equation}
\begin{equation}\nonumber
\begin{aligned}
{}&\Big|\frac{1}{4\pi^2}\int_{\mathbb{R}^2}\Big(\int_{\mathbb{R}^2}\frac{x_2-y_2}{|x-y|^2}W_{\varepsilon,Y}(y)
u(y)\,dy\cdot\int_{\mathbb{R}^2}\frac{x_2-z_2}{|x-z|^2}W_{\varepsilon,Y}(z)
v(z)\,dz\Big)W_{\varepsilon,Y}^2(x)\,dx\Big|\\
\leq{}& C\int_{\mathbb{R}^2}\Big(\int_{\mathbb{R}^2}\frac{1}{|x-y|^2}W_{\varepsilon,Y}^2(y)
\,dy\Big)^{\frac{1}{2}}\Big(\int_{\mathbb{R}^2}|u|^2\,dy\Big)^{\frac{1}{2}}\Big(\int_{\mathbb{R}^2}\frac{1}{|x-z|^2}
W_{\varepsilon,Y}^2(z)\,dz\Big)^{\frac{1}{2}}\Big(\int_{\mathbb{R}^2}|v|^2\,dz\Big)^{\frac{1}{2}}W_{\varepsilon,Y}^2(x)\,dx\\
\leq{}& C||u||_{\varepsilon}||v||_{\varepsilon}\int_{\mathbb{R}^2}W_{\varepsilon,Y}^2(x)\,dx
\leq C\varepsilon^2||u||_{\varepsilon}||v||_{\varepsilon},
\end{aligned}
\end{equation}
and
\begin{equation}\nonumber
\begin{aligned}
{}&\Big|\frac{1}{2\pi^2}\int_{\mathbb{R}^2}\Big(\int_{\mathbb{R}^2}\frac{x_2-y_2}{|x-y|^2}W_{\varepsilon,Y}^2(y)\,dy
\cdot\int_{\mathbb{R}^2}\frac{x_2-z_2}{|x-z|^2}W_{\varepsilon,Y}(z)u(z)\,dz\Big)W_{\varepsilon,Y}(x)v(x)\,dx\Big|\\
\leq {}&C\varepsilon\int_{\mathbb{R}^2}\Big(\int_{\mathbb{R}^2}\frac{1}{|x-z|}W_{\varepsilon,Y}(z)|u(z)|\,dz\Big)W_{\varepsilon,Y}(x)|v(x)|\,dx\\
\leq {}&C\varepsilon\int_{\mathbb{R}^2}\Big(\int_{\mathbb{R}^2}\frac{1}{|x-z|^2}W_{\varepsilon,Y}^2(z)\,dz\Big)^{\frac{1}{2}}
\Big(\int_{\mathbb{R}^2}|u|^2\,dy\Big)^{\frac{1}{2}}W_{\varepsilon,Y}(x)|v(x)|\,dx\\
\leq{}& C\varepsilon||u||_{\varepsilon}\int_{\mathbb{R}^2}W_{\varepsilon,Y}(x)|v(x)|\,dx
\leq C\varepsilon||u||_{\varepsilon}\Big(\int_{\mathbb{R}^2}W_{\varepsilon,Y}^2(x)\,dx\Big)^{\frac{1}{2}}\Big(\int_{\mathbb{R}^2}|v|^2\,dx\Big)^{\frac{1}{2}}
\leq C\varepsilon^2||u||_{\varepsilon}||v||_{\varepsilon}.
\end{aligned}
\end{equation}
Similarly, we have estimates of the rest terms of $|\langle L_{2,\varepsilon}u,v\rangle_\varepsilon|$. Hence, $L_{\varepsilon}$ is bounded.
\end{proof}



Now we are in a proposition to prove Proposition \ref{prop1}.
\noindent
\begin{proof}[\textbf{Proof of Proposition \ref{prop1}.}]
We will use the contraction mapping theorem to prove it.
As we all know, for fixed $Y\in D_{\delta}$, finding a critical point for $\mathcal{J}_{\varepsilon}(Y,\varphi)$ is equivalent to solving
\begin{equation}\label{equ1}
l_{\varepsilon}+L_{\varepsilon}\varphi+R'_{\varepsilon}(\varphi)=0.
\end{equation}
By Lemma \ref{3.4}, $L_{\varepsilon}$ is invertible in $E_{\varepsilon,Y}$. Thus solving \eqref{equ1} is equivalent
to finding a fixed point of
\begin{equation}
\varphi=A\varphi:=-L_{\varepsilon}^{-1}(l_{\varepsilon}+R'_{\varepsilon}(\varphi)).
\end{equation}

We set
\begin{equation}\nonumber
S_{\varepsilon}:=\Big\{\varphi\in E_{\varepsilon,Y}:||\varphi||_{\varepsilon}\leq\varepsilon\Big(\varepsilon^{\theta-\tau}+\sum_{i=1}^k|V(y^i)-V(a^i)|^{1-\tau}\Big)\Big\}
\end{equation}
where $\tau>0$ is a fixed small constant. Firstly, $A$ maps $S_{\varepsilon}$ to $S_{\varepsilon}$. In fact, for $\varphi\in S_{\varepsilon}$, by Lemma \ref{3.1} and \ref{3.2}, we obtain
\begin{equation}\nonumber
\begin{aligned}
||A\varphi||_\varepsilon\leq{}& C(||l_{\varepsilon}||+||R'_{\varepsilon}(\varphi)||)\\
\leq {}&C\Big(\varepsilon^{1+\theta}+\varepsilon\sum_{i=1}^k|V(y^i)-V(a^i)|+\varepsilon^{-\min\{1,p-2\}}||\varphi||_{\varepsilon}^{\min\{2,p-1\}}\\
{}&+\varepsilon||\varphi||_{\varepsilon}^{2}+
||\varphi||_{\varepsilon}^{3}+\varepsilon^{-1}||\varphi||_{\varepsilon}^{4}+\varepsilon^{-2}||\varphi||_{\varepsilon}^{5}\Big)\\
\leq{}& C\varepsilon\Big(\varepsilon^{\theta-\tau}+\sum_{i=1}^k|V(y^i)-V(a^i)|^{1-\tau}\Big).
\end{aligned}
\end{equation}
Secondly, $A$ is a contraction map. In fact, for any $\varphi_1,\varphi_2\in S_{\varepsilon}$, we have
\begin{equation}\nonumber
\begin{aligned}
||A\varphi_1-A\varphi_2||_{\varepsilon}&\leq C||R'_{\varepsilon}(\varphi_1)-R'_{\varepsilon}(\varphi_2)||\\
&= C||R''_{\varepsilon}(\sigma\varphi_1+(1-\sigma)\varphi_2)||||\varphi_1-\varphi_2||_{\varepsilon}\\
&\leq \frac{1}{2}||\varphi_1-\varphi_2||_{\varepsilon}.
\end{aligned}
\end{equation}
Therefore, by the contraction mapping theorem, we can conclude that for fixed $Y\in D_{\delta}$, $\mathcal{J}_{\varepsilon}(Y,\varphi)$ has a unique critical point. Thus
\begin{equation}\nonumber
\begin{aligned}
{}&||\varphi_{\varepsilon,Y}||_{\varepsilon}=||A\varphi_{\varepsilon,Y}||_{\varepsilon}
\leq C(||l_{\varepsilon}||+||R'_{\varepsilon,Y}(\varphi)||)\\
\leq{}& C||l_{\varepsilon}||
+C\Big((\varepsilon^{-1}||\varphi_{\varepsilon,Y}||_{\varepsilon})^{\min\{1,p-2\}}+\varepsilon||\varphi_{\varepsilon,Y}||_{\varepsilon}+
||\varphi_{\varepsilon,Y}||_{\varepsilon}^{2}+\varepsilon^{-1}||\varphi_{\varepsilon,Y}||_{\varepsilon}^{3}+\varepsilon^{-2}
||\varphi_{\varepsilon,Y}||_{\varepsilon}^{4}\Big)||\varphi_{\varepsilon,Y}||_{\varepsilon},\\
\end{aligned}
\end{equation}
which yields
\begin{equation}\nonumber
||\varphi_{\varepsilon,Y}||_{\varepsilon}\leq C||l_{\varepsilon}||\leq C\varepsilon
\Big(\varepsilon^\theta+\sum_{i=1}^k|V(y^i)-V(a^i)|\Big).
\end{equation}
This completes the proof of Proposition \ref{prop1}.
\end{proof}

\subsection{Proof of the existence of  k-peak solutions}\label{sec3}
\noindent

In this subsection, we mainly prove Theorem \ref{thm1}.


\begin{proof}[\textbf{Proof of Theorem \ref{thm1}.}] Let $\varepsilon_0$ and $\delta_0$ be defined as in Proposition \ref{prop1} and let $\varepsilon<\varepsilon_0$. Fix $0<\delta<\delta_0$. Let $Y\mapsto\varphi_{\varepsilon,Y}$ for $Y\in D_\delta$ be the map obtained in Proposition \ref{prop1}. We will find a critical point for the function $j_\varepsilon(Y):=\mathcal{J}_{\varepsilon}(Y,\varphi_{\varepsilon,Y})$.
By \eqref{eq-J} and Lemmas \ref{3.1} to \ref{3.3}, we have
\begin{align*}
j_\varepsilon(Y)={}&I_{\varepsilon}(W_{\varepsilon,Y})+\langle l_{\varepsilon},\varphi_{\varepsilon,Y}\rangle_\varepsilon+\frac{1}{2}
\langle L_{\varepsilon}\varphi_{\varepsilon,Y},\varphi_{\varepsilon,Y}\rangle_\varepsilon+R_{\varepsilon}(\varphi_{\varepsilon,Y})\\
={}&I_{\varepsilon}(W_{\varepsilon,Y})+O(||l_{\varepsilon}||||\varphi_{\varepsilon,Y}||_\varepsilon+
||\varphi_{\varepsilon,Y}||_\varepsilon^2)\\
={}&(\frac{1}{2}-\frac{1}{p})\varepsilon^2\sum_{i=1}^k\int_{\mathbb{R}^2}(U^i(x))^p\,dx
+\frac{1}{2}\varepsilon^2\sum_{i=1}^k\int_{\mathbb{R}^2}(V(y^i)-V(a^i))(U^i(x))^2\,dx\\
{}&+O(\varepsilon^{2+\theta})
+O\Big(\varepsilon^2\Big(\varepsilon^{\theta}+\sum_{i=1}^k|V(y^i)-V(a^i)|\Big)^2\Big)\\
=:{}&C_1\varepsilon^2+\varepsilon^2\sum_{i=1}^kC_{2,i}(V(y^i)-V(a^i))+O(\varepsilon^{2+\theta})+
O\Big(\varepsilon^2\Big(\varepsilon^{\theta}+\sum_{i=1}^k|V(y^i)-V(a^i)|\Big)^2\Big)
\end{align*}
where $C_1=(\frac{1}{2}-\frac{1}{p})\sum_{i=1}^k\int_{\mathbb{R}^2}(U^i(x))^p\,dx$ and $C_{2,i}=\frac{1}{2}\int_{\mathbb{R}^2}(U^i(x))^2\,dx$.

Consider the following minimizing problem
\begin{equation}
j_{\varepsilon}(Y_{\varepsilon})\equiv\min_{Y\in D_{\delta}}j_{\varepsilon}(Y).
\end{equation}
Applying a comparison argument, we claim that $Y_{\varepsilon}$ is an interior point of $D_{\delta}$ and hence $Y_{\varepsilon}$ is a critical point of $j_{\varepsilon}(Y)$ for $\varepsilon$ is sufficiently small.

To prove this, choose $e^i\in\mathbb{R}^2$ for $i=1,\cdots,k,$ with $|e^i|=1,\ e^i\neq e^j$ for $i\neq j$ and $\eta>1$. Let $z_\varepsilon^i=a^i+\varepsilon^\eta e^i$ satisfying $Z_\varepsilon=(z_\varepsilon^1,\cdots,z_\varepsilon^k)\in D_\delta,$ where $\eta$ is a sufficiently large constant. Thus applying the H\"{o}lder continuity of $V$, we have
\begin{equation}\nonumber
\begin{aligned}
j_{\varepsilon}(Z_{\varepsilon})&=C_1\varepsilon^2+O(\varepsilon^{2+\eta\theta}+\varepsilon^
{2+\theta}+\varepsilon^2(\varepsilon^\theta+\varepsilon^{\eta\theta})^2)\\
&=C_1\varepsilon^2+O(\varepsilon^{2+\theta}).
\end{aligned}
\end{equation}
On the other hand, by using $j_{\varepsilon}(Y_{\varepsilon})\leq j_{\varepsilon}(Z_{\varepsilon})$, we obtain
\begin{equation}\nonumber
\varepsilon^2\sum_{i=1}^kC_{2,i}(V(y_\varepsilon^i)-V(a^i))+
O\Big(\varepsilon^2\Big(\varepsilon^{\theta}+\sum_{i=1}^k|V(y_\varepsilon^i)-V(a^i)|\Big)^2\Big)\leq O(\varepsilon^{2+\theta}).
\end{equation}
If $Y_\varepsilon\in\partial D_\delta$, by the assumption $(V_2)$ we get
\begin{equation}\nonumber
V(y_\varepsilon^i)-V(a^i)\geq c_i>0,\ j=1,\cdots,k
\end{equation}
for some constants $c_i>0$. Then
\begin{equation}\nonumber
0<C_{2,i}\leq0.
\end{equation}
This leads to a contradiction. Therefore $Y_{\varepsilon}$ is an interior point of $D_{\delta}$.
\end{proof}

\section{Local uniqueness of the k-peak solutions}\label{sec3-1}
First we give the concrete form of concentrated solutions  for \eqref{eq1}.

\begin{prop}\label{prop2}
Let $\{u_{\varepsilon}(x)\}_{\epsilon>0}$ be a family of positive solutions of \eqref{eq1}
concentrating at different points $\{ a^1,\ldots,a^k\} \subset\mathbb{R}^2$ with $k\geq 1.$ Then, for $\varepsilon>0$ sufficiently small, $u_\varepsilon$ is of the form
\begin{equation}\label{101}
u_\varepsilon(x)=\sum_{i=1}^kU^i(\frac{x-y_\varepsilon^i}{\varepsilon})+\varphi_\varepsilon(x)
\end{equation}
with
\begin{equation}\label{102}
|y_\varepsilon^i-a^i|=o(1)\quad \text{and}\quad ||\varphi_\varepsilon||_{\varepsilon}=o(\varepsilon),
\end{equation}
and
 \begin{equation}\label{103-add}
\Big(\varphi_\varepsilon(x),\frac{ \partial U_{i}\big(\frac{x-y^{i}_{\epsilon}}{\epsilon}\big)}{\partial x^{i}}\Big)_{\varepsilon}=0,\,\,\,\,i=1,2.
\end{equation}
\end{prop}

\noindent
\begin{proof}
For each $1\leq i\leq k$, let $u^i_{\varepsilon}(x)=u_\varepsilon(\varepsilon x+y_\varepsilon^i)$. Then $u^i_{\varepsilon}$ is a uniformly bounded sequence in $H^1(\mathbb{R}^2)$ with respect to $\varepsilon$ and satisfies
\begin{equation}\label{103}\nonumber
-\Delta u^i_{\varepsilon}+V(\varepsilon x+y_\varepsilon^i)u^i_{\varepsilon}+\varepsilon^2(A_0(u^i_{\varepsilon})+A_1^2(u^i_{\varepsilon})+A_2^2(
u^i_{\varepsilon}))u^i_{\varepsilon}=(u^i_{\varepsilon})^{p-1},\ x\in\mathbb{R}^2.
\end{equation}
Suppose that $\psi(x)$ is an arbitrarily fixed function in $H^1(\mathbb{R}^2)$. Then
\begin{equation}\label{104}
\int_{\mathbb{R}^2}\left(\nabla u^i_{\varepsilon}\nabla\psi+V(\varepsilon x+y_\varepsilon^i)u^i_{\varepsilon}\psi+\varepsilon^2(A_0(u^i_{\varepsilon})+A_1^2(u^i_{\varepsilon})+A_2^2(
u^i_{\varepsilon}))u^i_{\varepsilon}\psi\right)\,dx=\int_{\mathbb{R}^2}(u^i_{\varepsilon})^{p-1}\psi\,dx.
\end{equation}
By Sobolev embedding, similar to \eqref{equa3}, we can deduce
\begin{equation}\nonumber
\begin{aligned}
A_1^2(u_\varepsilon^i)\leq{}&C \Big(\int_{B_{1}(x)}\frac{1}{|x-y|}(u^i_\varepsilon(y))^2\,dy\Big)^2+C\Big(\int_{\mathbb{R}^2\backslash B_1(x)}\frac{1}{|x-y|}(u^i_\varepsilon(y))^2\,dy\Big)^2\\
\leq{}&C\Big(\int_{B_{1}(x)}\frac{1}{|x-y|^{\frac{3}{2}}}\,dy\Big)^{\frac{4}{3}}
\Big(\int_{B_1(x)}(u^i_\varepsilon(y))^6\,dy\Big)^{\frac{2}{3}}+C\Big(\int_{\mathbb{R}^2}(u^i_\varepsilon(y))^2\,dy\Big)^2\\
\leq{}&C||u^i_\varepsilon||^4_{H^1(\mathbb{R}^2)}\leq C.\\
\end{aligned}
\end{equation}
Analogously, we get $A_0(u^i_{\varepsilon})\leq C$ and $A_2^2(
u^i_{\varepsilon})\leq C$.
Passing to a subsequence if necessary, there exists $w^i\in H^1(\mathbb{R}^2)$ such that $u^i_{\varepsilon}\rightharpoonup w^i$ in $H^1(\mathbb{R}^2)$. Letting $\varepsilon\rightarrow0$ in \eqref{104}, we get
\begin{equation}\label{105}\nonumber
\int_{\mathbb{R}^2}\left(\nabla w^i\nabla\psi+V(a^i)w^i\psi\right)\,dx=\int_{\mathbb{R}^2}(w^i)^{p-1}\psi\,dx,
\end{equation}
that is, $w^i$ is a weak solution of
\begin{equation}\label{106}\nonumber
-\Delta u+V(a^i)u=u^{p-1}.
\end{equation}
Note that $x=0$ is a maximum point of $w^i$. Hence $w^i(x)=U^i(|x|)$ must be the unique positive radial solution of \eqref{rr}. By the same concentrating compactness arguments as in \cite{ref27}, we can write $u_\varepsilon(x)$ uniquely as \eqref{101} with $y^i_\varepsilon$ and $\varphi_\varepsilon$ satisfying \eqref{102}.
\end{proof}


%
%

Now we obtain the more precise estimates of $|y_\varepsilon^i-a^i|$ and $||\varphi_\varepsilon||_{\varepsilon}.$
\begin{prop}\label{prop3}
Suppose that $p>2$ and $V(x)$ satisfies $(V_1),\ (V_2)$ and $(V_3)$. Let $u_{\epsilon}(x)$  be the solution of \eqref{eq1}
concentrating at $k(k\geq 1)$ different local minimum points $\{ a^1,\ldots,a^k\} \subset\mathbb{R}^2$ of $V(x).$
 Suppose that
\begin{equation}\label{107}
u_\varepsilon(x)=\sum_{i=1}^kU^i(\frac{x-y_\varepsilon^i}{\varepsilon})+\varphi_\varepsilon(x).
\end{equation}
Then
\begin{equation}\label{108}
|y_\varepsilon^i-a^i|=o(\varepsilon)
\end{equation}
for each $i=1,\cdots,k,$ and
\begin{equation}\label{109}
||\varphi_\varepsilon||_{\varepsilon}=O(\varepsilon^{1+m}).
\end{equation}
\end{prop}

\noindent
\begin{proof}
Choose $i$ such that
\begin{equation}\label{110}\nonumber
|y_\varepsilon^{i}-a^{i}|=\max\{|y_\varepsilon^j-a^j|:1\leq j\leq k\}.
\end{equation}
Since $V(x)$ satisfies $(V_1),\ (V_2)$ and $(V_3)$, similar to the proof of Lemma \ref{3.1}, we have
\begin{equation}\label{111}\nonumber
||l_\varepsilon||\leq C\varepsilon(\varepsilon^{m}+|y_{\varepsilon}^i-a^i|^m).
\end{equation}
As a result, following the proof of Proposition \ref{prop1}, there holds
\begin{equation}\label{112}
||\varphi_\varepsilon||_\varepsilon\leq C||l_\varepsilon||\leq C\varepsilon(\varepsilon^{m}+|y_{\varepsilon}^i-a^i|^m).
\end{equation}
Next it remains to estimate $|y_{\varepsilon}^i-a^i|$. Applying the Pohozaev identity \eqref{P} to $u=u_\varepsilon$ with $\Omega=B_d(y^i_\varepsilon)$ for some small constant $d>0$, we have
\begin{equation}\label{113}
\begin{aligned}
\int_{B_d(y^i_\varepsilon)}\frac{\partial V}{\partial x_k}u_\varepsilon^2\,dx
={}&\int_{\partial B_d(y^i_\varepsilon)}(\varepsilon^2|\nabla u_\varepsilon|^2+V(x)u_\varepsilon^2)\nu_k\,d\sigma-2\int_{\partial B_d(y^i_\varepsilon)}\varepsilon^2\frac{\partial u_\varepsilon}{\partial\nu}\frac{\partial u_\varepsilon}{\partial x_k}\,d\sigma\\
{}&-\frac{2}{p}\int_{\partial B_d(y^i_\varepsilon)}|u_\varepsilon|^p\nu_k\,d\sigma-2\int_{\partial B_d(y^i_\varepsilon)}(A_1\nu_1
+A_2\nu_2)A_ku_\varepsilon^2\,d\sigma\\
{}&+\int_{\partial B_d(y^i_\varepsilon)}
(A_0+A_1^2+A_2^2)u_\varepsilon^2\nu_k\,d\sigma.
\end{aligned}
\end{equation}
By the properties of $U^j$, we infer that
\begin{equation}\label{114}
\max_{1\leq j\leq k}(U^j(x)+|\nabla U^j(x)|)=O(e^{-\sigma |x|})
\end{equation}
for some $\sigma>0$. Then
\begin{equation}\label{115}\nonumber
\begin{aligned}
\int_{\partial B_d(y^i_\varepsilon)}\varepsilon^2|\nabla u_\varepsilon|^2\nu_k\,d\sigma
&\leq \varepsilon^2\int_{\partial B_d(y^i_\varepsilon)}|\nabla u_\varepsilon|^2\,d\sigma\\
&\leq C\varepsilon^2\int_{\partial B_d(y^i_\varepsilon)}\Big(\sum_{j=1}^k|\nabla U^j_{\varepsilon,y^j_\varepsilon}|^2+|\nabla\varphi_\varepsilon|^2\Big)\,d\sigma\\
&\leq C(e^{-\frac{\tau}{\varepsilon}}+||\varphi_\varepsilon||^2_\varepsilon)
\end{aligned}
\end{equation}
for some $\tau>0$, which implies that for any $\gamma>0$, there holds
\begin{equation}\label{116}
\int_{\partial B_d(y^i_\varepsilon)}\varepsilon^2|\nabla u_\varepsilon|^2\nu_k\,d\sigma=O(\varepsilon^\gamma+
||\varphi_\varepsilon||^2_\varepsilon).
\end{equation}
Similarly, we have
\begin{equation}\label{118}
\begin{aligned}
\int_{\partial B_d(y^i_\varepsilon)}V(x)u_\varepsilon^2\nu_k\,d\sigma
&=O\Big(\int_{\partial B_d(y^i_\varepsilon)}u_\varepsilon^2\,d\sigma\Big)\\
&=O\Big(\int_{\partial B_d(y^i_\varepsilon)}\Big(\sum_{j=1}^k|U^j_{\varepsilon,y^j_\varepsilon}|^2+|\varphi_\varepsilon|^2\Big)\,d\sigma\Big)\\
&=O(\varepsilon^\gamma+||\varphi_\varepsilon||^2_\varepsilon)
\end{aligned}
\end{equation}
and
\begin{equation}\label{117}
\int_{\partial B_d(y^i_\varepsilon)}\varepsilon^2\frac{\partial u_\varepsilon}{\partial\nu}\frac{\partial u_\varepsilon}{\partial x_k}\,d\sigma=O\Big(\varepsilon^2\int_{\partial B_d(y^i_\varepsilon)}|\nabla u_\varepsilon|^2\,d\sigma\Big)=O(\varepsilon^\gamma+
||\varphi_\varepsilon||^2_\varepsilon).
\end{equation}
By Proposition \ref{prop2}, we know that $||\varphi_\varepsilon||_{\varepsilon}=o(\varepsilon)$. Then
\begin{equation}\label{119}
\begin{aligned}
\int_{\partial B_d(y^i_\varepsilon)}|u_\varepsilon|^p\nu_k\,d\sigma
&=O\Big(\int_{\partial B_d(y^i_\varepsilon)}\Big(\sum_{j=1}^k|U^j_{\varepsilon,y^j_\varepsilon}|^p+|\varphi_\varepsilon|^p\Big)\,d\sigma\Big)\\
&=O(\varepsilon^\gamma+\varepsilon^{2-p}||\varphi_\varepsilon||^p_{\varepsilon})\\
&=O(\varepsilon^\gamma+||\varphi_\varepsilon||^2_\varepsilon).
\end{aligned}
\end{equation}
According to \eqref{equa3}, we can verify that for $k=1,2$,
\begin{equation}\label{120}
\begin{aligned}
A_k(u_\varepsilon(x))&=O\Big(\int_{\mathbb{R}^2}\frac{1}{|x-y|}u_\varepsilon^2(y)\,dy\Big)\\
&=O\Big(\sum_{j=1}^k\int_{\mathbb{R}^2}\frac{1}{|x-y|}(U^j_{\varepsilon,y^j_\varepsilon}(y))^2\,dy\Big)
+O\Big(\int_{\mathbb{R}^2}\frac{1}{|x-y|}\varphi_\varepsilon^2(y)\,dy\Big)\\
&=O(\varepsilon+\varepsilon^{-1}||\varphi_\varepsilon||^2_\varepsilon),
\end{aligned}
\end{equation}
and
\begin{equation}\label{121}
\begin{aligned}
A_0(u_\varepsilon(x))&=O\Big(\int_{\mathbb{R}^2}\frac{1}{|x-y|}u_\varepsilon^2(y)
\Big(\int_{\mathbb{R}^2}\frac{1}{|y-z|}u_\varepsilon^2(z)\,dz\Big)\,dy\Big)\\
&=O\Big((\varepsilon+\varepsilon^{-1}||\varphi_\varepsilon||^2_\varepsilon)
\int_{\mathbb{R}^2}\frac{1}{|x-y|}u_\varepsilon^2(y)\,dy\Big)\\
&=O(\varepsilon^2+\varepsilon^{-2}||\varphi_\varepsilon||^4_\varepsilon).
\end{aligned}
\end{equation}
As a result,
\begin{equation}\label{122}\nonumber
\quad\int_{\partial B_d(y^i_\varepsilon)}(A_1\nu_1
+A_2\nu_2)A_ku_\varepsilon^2\,d\sigma
=O\Big((\varepsilon^2+\varepsilon^{-2}||\varphi_\varepsilon||^4_\varepsilon)
\int_{\partial B_d(y^i_\varepsilon)}u_\varepsilon^2\,d\sigma\Big)
=O(\varepsilon^\gamma+||\varphi_\varepsilon||^2_\varepsilon),
\end{equation}
and similarly,
\begin{equation}\label{123}\nonumber
\int_{\partial B_d(y^i_\varepsilon)}
(A_0+A_1^2+A_2^2)u_\varepsilon^2\nu_k\,d\sigma=O(\varepsilon^\gamma+||\varphi_\varepsilon||^2_\varepsilon).
\end{equation}
Therefore combining above estimates we obtain
\begin{equation}\label{124}
\text{RHS of \eqref{113}}=O(\varepsilon^\gamma+||\varphi_\varepsilon||^2_\varepsilon).
\end{equation}
On the other hand, we have
\begin{equation}\label{125}
\text{LHS of \eqref{113}}=
\int_{B_d(y^i_\varepsilon)}\frac{\partial V}{\partial x_k}(U^i_{\varepsilon,y^i_\varepsilon})^2\,dx
+2\int_{B_d(y^i_\varepsilon)}\frac{\partial V}{\partial x_k}U^i_{\varepsilon,y^i_\varepsilon}\varphi_\varepsilon\,dx+O(\varepsilon^\gamma+||\varphi_\varepsilon||^2_\varepsilon).
\end{equation}
It follows from the assumption $(V_3)$ that
\begin{equation}\label{126}
\begin{aligned}
\int_{B_d(y^i_\varepsilon)}\frac{\partial V}{\partial x_k}U^i_{\varepsilon,y^i_\varepsilon}\varphi_\varepsilon\,dx
&\leq C\Big(\int_{B_d(y^i_\varepsilon)}(\frac{\partial V}{\partial x_k}U^i_{\varepsilon,y^i_\varepsilon})^2\,dx\Big)^\frac{1}{2}||\varphi_\varepsilon||_\varepsilon\\
&\leq C\Big(\varepsilon^2\int_{B_{\frac{d}{\varepsilon}}(0)}(\frac{\partial V}{\partial x_k}(\varepsilon x+y^i_\varepsilon)U^i)^2\,dx\Big)^\frac{1}{2}||\varphi_\varepsilon||_\varepsilon\\
&\leq C\Big(\varepsilon^2\int_{B_{\frac{d}{\varepsilon}}(0)}(|\varepsilon x+y^i_\varepsilon-a^i|^{m-1}
U^i)^2\,dx\Big)^\frac{1}{2}||\varphi_\varepsilon||_\varepsilon\\
&\leq C\varepsilon(\varepsilon^{m-1}+|y^i_\varepsilon-a^i|^{m-1})||\varphi_\varepsilon||_\varepsilon.
\end{aligned}
\end{equation}
Also, we have
\begin{equation}\label{127}
\begin{aligned}
\int_{B_d(y^i_\varepsilon)}\frac{\partial V}{\partial x_k}(U^i_{\varepsilon,y^i_\varepsilon})^2\,dx
={}&mb_{k,i}\int_{B_d(y^i_\varepsilon)}|x_k-a^i_k|^{m-2}(x_k-a^i_k)(U^i_{\varepsilon,y^i_\varepsilon})^2\,dx\\
{}&+O\Big(\int_{B_d(y^i_\varepsilon)}|x-a^i|^m(U^i_{\varepsilon,y^i_\varepsilon})^2\,dx\Big)\\
={}&mb_{k,i}\varepsilon^2\int_{B_{\frac{d}{\varepsilon}}(0)}|\varepsilon x_k+y^i_{\varepsilon,k}-a^i_k|^{m-2}(\varepsilon x_k+y^i_{\varepsilon,k}-a^i_k)(U^i)^2\,dx\\
{}&+O(\varepsilon^2(\varepsilon^{m}+|y^i_\varepsilon-a^i|^{m})).
\end{aligned}
\end{equation}
Combining \eqref{124}-\eqref{127} and \eqref{112} yields
\begin{equation}\label{128}
mb_{k,i}\int_{B_{\frac{d}{\varepsilon}}(0)}|\varepsilon x_k+y^i_{\varepsilon,k}-a^i_k|^{m-2}(\varepsilon x_k+y^i_{\varepsilon,k}-a^i_k)(U^i)^2\,dx=O(\varepsilon^{m}+|y^i_\varepsilon-a^i|^{m})
\end{equation}
by choosing $\gamma>0$ sufficiently large. Note that
\begin{equation}\label{129}\nonumber
||a+b|^m-|a|^m-m|a|^{m-2}ab|\leq C(|a|^{m-m^\ast}|b|^{m^\ast}+|b|^m),
\end{equation}
where $a,b\in\mathbb{R},m>1,m^\ast=\min\{m,2\}$ and the constant $C$ is independent of $a$ and $b$.
Taking $a=\varepsilon x_k+y^i_{\varepsilon,k}-a^i_k$ and $b=-\varepsilon x_k$, we obtain
\begin{equation}\label{130}\nonumber
\begin{aligned}
{}& m|\varepsilon x_k+y^i_{\varepsilon,k}-a^i_k|^{m-2}(\varepsilon x_k+y^i_{\varepsilon,k}-a^i_k)(y^i_{\varepsilon,k}-a^i_k)\\
={}&m|\varepsilon x_k+y^i_{\varepsilon,k}-a^i_k|^{m}-m|\varepsilon x_k+y^i_{\varepsilon,k}-a^i_k|^{m-2}(\varepsilon x_k+y^i_{\varepsilon,k}-a^i_k)\varepsilon x_k\\
\geq{}& |y^i_{\varepsilon,k}-a^i_k|^m+(m-1)|\varepsilon x_k+y^i_{\varepsilon,k}-a^i_k|^{m}-C(
|\varepsilon x_k+y^i_{\varepsilon,k}-a^i_k|^{m-m^\ast}|\varepsilon x_k|^{m^\ast}+|\varepsilon x_k|^m)\\
\geq{}& |y^i_{\varepsilon,k}-a^i_k|^m-C(
|y^i_{\varepsilon,k}-a^i_k|^{m-m^\ast}|\varepsilon x_k|^{m^\ast}+|\varepsilon x_k|^m).
\end{aligned}
\end{equation}
So we get
\begin{equation}\label{131}
\begin{aligned}
{}&|y^i_{\varepsilon,k}-a^i_k|\Big|\int_{B_{\frac{d}{\varepsilon}}(0)}|\varepsilon x_k+y^i_{\varepsilon,k}-a^i_k|^{m-2}(\varepsilon x_k+y^i_{\varepsilon,k}-a^i_k)(U^i)^2\,dx\Big|\\
\geq{}& \frac{1}{m}|y^i_{\varepsilon,k}-a^i_k|^m\int_{B_{\frac{d}{\varepsilon}}(0)}(U^i)^2\,dx-\frac{C}{m}
\int_{B_{\frac{d}{\varepsilon}}(0)}(|y^i_{\varepsilon,k}-a^i_k|^{m-m^\ast}|\varepsilon x_k|^{m^\ast}+|\varepsilon x_k|^m)(U^i)^2\,dx\\
\geq{}& \frac{1}{m}|y^i_{\varepsilon,k}-a^i_k|^m\int_{B_{\frac{d}{\varepsilon}}(0)}(U^i)^2\,dx-\frac{C}{m}
(\varepsilon^m+\varepsilon^{m^\ast}|y^i_{\varepsilon,k}-a^i_k|^{m-m^\ast}).
\end{aligned}
\end{equation}
Choose $k_0\in\{1,2\}$ such that $|y^i_{\varepsilon,k_0}-a^i_{k_0}|\geq\frac{|y^i_{\varepsilon,k}-a^i_k|}{\sqrt{2}}$. Using \eqref{128} and \eqref{131} and applying $\varepsilon$-Young inequality, we get
\begin{equation}\label{132}\nonumber
\begin{aligned}
|y^i_{\varepsilon}-a^i|^m&=|y^i_{\varepsilon}-a^i|O(\varepsilon^{m}+|y^i_\varepsilon-a^i|^{m})
+O(\varepsilon^m+\varepsilon^{m^\ast}|y^i_{\varepsilon}-a^i|^{m-m^\ast})\\
&=|y^i_{\varepsilon}-a^i|O(\varepsilon^{m}+|y^i_\varepsilon-a^i|^{m})
+O(\varepsilon^m)+\frac{1}{2}|y^i_{\varepsilon}-a^i|^m.
\end{aligned}
\end{equation}
Thus $|y^i_{\varepsilon}-a^i|=O(\varepsilon)$.
Taking a subsequence necessarily, we may assume that $\frac{y^i_{\varepsilon}-a^i}{\varepsilon}\to t=(t_1,t_2)\in\mathbb{R}^2$ as $\varepsilon\to0$. By \eqref{128} we have
\begin{equation}\label{134}\nonumber
\int_{B_{\frac{d}{\varepsilon}}(0)}\Big|x_k+\frac{y^i_{\varepsilon,k}-a^i_k}{\varepsilon}\Big|^{m-2}\Big( x_k+\frac{y^i_{\varepsilon,k}-a^i_k}{\varepsilon}\Big)(U^i)^2\,dx=O(\varepsilon)
\end{equation}
for $k=1,2$, which gives
\begin{equation}\label{135}\nonumber
\int_{\mathbb{R}^2}|x_k+t_k|^{m-2}(x_k+t_k)(U^i)^2\,dx=0.
\end{equation}
Since $U^i$ is radially symmetric decreasing, we get $t=0$. Thus $|y^i_{\varepsilon}-a^i|=o(\varepsilon)$ and from \eqref{112} we obtain $||\varphi_\varepsilon||_\varepsilon=O(\varepsilon^{m+1})$.
\end{proof}

Now suppose that $u_\varepsilon^{(j)}(x)=\sum_{i=1}^kU^i(\frac{x-y_\varepsilon^{i(j)}}{\varepsilon})+\varphi_\varepsilon^{(j)}(x),\ j=1,2,$ are two k-peak solutions of \eqref{eq1} defined as in the Definition \ref{def1} concentrating around $a^i$ for $i=1,\cdots,k$. Define
\begin{equation}\label{136}\nonumber
\xi_\varepsilon=\frac{u_\varepsilon^{(1)}-u_\varepsilon^{(2)}}{||u_\varepsilon^{(1)}-u_\varepsilon^{(2)}||_{L^\infty(\mathbb{R}^2)}}.
\end{equation}
Then $\xi_\varepsilon$ satisfies $||\xi_\varepsilon||_{L^\infty(\mathbb{R}^2)}=1$ and
\begin{equation}\label{137}
-\varepsilon^2\Delta\xi_\varepsilon+V(x)\xi_\varepsilon+(A_0(u_\varepsilon^{(1)})+ A_1^2(u_\varepsilon^{(1)})+A_2^2(u_\varepsilon^{(1)}))\xi_\varepsilon=C_\varepsilon(x)\xi_\varepsilon+E(x),
\end{equation}
where
\begin{equation}\label{138}\nonumber
C_\varepsilon(x)=(p-1)\int_0^1(tu_\varepsilon^{(1)}(x)+(1-t)u_\varepsilon^{(2)}(x))^{p-2}\,dt
\end{equation}
and
\begin{align*}
E(x)={}&\frac{u_\varepsilon^{(2)}(x)}{8\pi^2}\int_{\mathbb{R}^2}\frac{x_1-y_1}{|x-y|^2}(u_\varepsilon^{(1)}(y)+u_\varepsilon^{(2)}(y))
\xi_\varepsilon(y)\Big(\int_{\mathbb{R}^2}\frac{y_1-z_1}{|y-z|^2}(u_\varepsilon^{(2)}(z))^2\,dz\Big)\,dy\\
{}&+\frac{u_\varepsilon^{(2)}(x)}{8\pi^2}\int_{\mathbb{R}^2}\frac{x_1-y_1}{|x-y|^2}(u_\varepsilon^{(1)}(y))^2
\Big(\int_{\mathbb{R}^2}\frac{y_1-z_1}{|y-z|^2}(u_\varepsilon^{(1)}(z)+u_\varepsilon^{(2)}(z))\xi_\varepsilon(z)\,dz\Big)\,dy\\
{}&-\frac{u_\varepsilon^{(2)}(x)}{16\pi^2}\int_{\mathbb{R}^2}\frac{x_1-y_1}{|x-y|^2}((u_\varepsilon^{(1)}(y))^2+(u_\varepsilon^{(2)}(y))^2)\,dy
\cdot\int_{\mathbb{R}^2}\frac{x_1-z_1}{|x-z|^2}(u_\varepsilon^{(1)}(z)+u_\varepsilon^{(2)}(z))\xi_\varepsilon(z)\,dz\\
{}&+\frac{u_\varepsilon^{(2)}(x)}{8\pi^2}\int_{\mathbb{R}^2}\frac{x_2-y_2}{|x-y|^2}(u_\varepsilon^{(1)}(y)+u_\varepsilon^{(2)}(y))
\xi_\varepsilon(y)\Big(\int_{\mathbb{R}^2}\frac{y_2-z_2}{|y-z|^2}(u_\varepsilon^{(2)}(z))^2\,dz\Big)\,dy\\
{}&+\frac{u_\varepsilon^{(2)}(x)}{8\pi^2}\int_{\mathbb{R}^2}\frac{x_2-y_2}{|x-y|^2}(u_\varepsilon^{(1)}(y))^2
\Big(\int_{\mathbb{R}^2}\frac{y_2-z_2}{|y-z|^2}(u_\varepsilon^{(1)}(z)+u_\varepsilon^{(2)}(z))\xi_\varepsilon(z)\,dz\Big)\,dy\\
{}&-\frac{u_\varepsilon^{(2)}(x)}{16\pi^2}\int_{\mathbb{R}^2}\frac{x_2-y_2}{|x-y|^2}((u_\varepsilon^{(1)}(y))^2+(u_\varepsilon^{(2)}(y))^2)\,dy
\cdot\int_{\mathbb{R}^2}\frac{x_2-z_2}{|x-z|^2}(u_\varepsilon^{(1)}(z)+u_\varepsilon^{(2)}(z))\xi_\varepsilon(z)\,dz\\
=:{}&E_{1}(x)+E_{2}(x)+E_{3}(x)+E_{4}(x)+E_{5}(x)+E_{6}(x).
\end{align*}
In the rest of this section, we aim to prove $||\xi_\varepsilon||_{L^\infty(\mathbb{R}^2)}=o(1)$, which is in contradiction with $||\xi_\varepsilon||_{L^\infty(\mathbb{R}^2)}=1$.

\begin{prop}\label{prop7}
There holds
\begin{equation}\label{146}
||\xi_\varepsilon||_\varepsilon=O(\varepsilon).
\end{equation}
\end{prop}

\noindent
\begin{proof}
From \eqref{137}, we have
\begin{equation}\label{147}\nonumber
||\xi_\varepsilon||_\varepsilon^2=-\int_{\mathbb{R}^2}(A_0(u_\varepsilon^{(1)})+ A_1^2(u_\varepsilon^{(1)})+A_2^2(u_\varepsilon^{(1)}))\xi^2_\varepsilon\,dx+\int_{\mathbb{R}^2}
C_\varepsilon(x)\xi^2_\varepsilon\,dx+\int_{\mathbb{R}^2}E(x)\xi_\varepsilon\,dx.
\end{equation}
Combining \eqref{109}, \eqref{120} and \eqref{121} yields
\begin{equation}\label{148}\nonumber
\int_{\mathbb{R}^2}(A_0(u_\varepsilon^{(1)})+ A_1^2(u_\varepsilon^{(1)})+A_2^2(u_\varepsilon^{(1)}))\xi^2_\varepsilon\,dx
=O(\varepsilon^2+\varepsilon^{-2}||\varphi_\varepsilon^{(1)}||
^4_\varepsilon)\int_{\mathbb{R}^2}\xi^2_\varepsilon\,dx=o(1)||\xi_\varepsilon||_\varepsilon^2.
\end{equation}
By calculating directly, we deduce
\begin{align*}
\int_{\mathbb{R}^2}C_\varepsilon(x)\xi^2_\varepsilon\,dx
&=O\Big(\sum_{i=1}^2\int_{\mathbb{R}^2}|u_\varepsilon^{(i)}|^{p-2}\xi^2_\varepsilon\,dx\Big)\\
&=O\Big(\sum_{i=1}^2\int_{\mathbb{R}^2}\Big(\sum_{j=1}^k|U^j_{\varepsilon,y_\varepsilon^{j(i)}}|^{p-2}
+|\varphi_\varepsilon^{(i)}|^{p-2}\Big)\xi^2_\varepsilon\,dx\Big)\\
&=O\Big(\sum_{i=1}^2\sum_{j=1}^k\varepsilon^2\int_{\mathbb{R}^2}|U^j|^{p-2}\,dx\Big)
+O\Big(\sum_{i=1}^2||\varphi_\varepsilon^{(i)}||^{p-2}_{L^p(\mathbb{R}^2)}||\xi_\varepsilon||^2_{L^p(\mathbb{R}^2)}\Big)\\
&=O(\varepsilon^2)+O\Big(\sum_{i=1}^2(\varepsilon^{2-p}||\varphi_\varepsilon^{(i)}||^{p}_\varepsilon)^
{\frac{p-2}{p}}(\varepsilon^{2-p}
||\xi_\varepsilon||^p_\varepsilon)^{\frac{2}{p}}\Big)\\
&=O(\varepsilon^2)+o(1)||\xi_\varepsilon||_\varepsilon^2.
\end{align*}
Using the estimate \eqref{120} and H\"{o}lder inequality, we can verify that
\begin{equation}\label{150}
\int_{\mathbb{R}^2}\frac{1}{|y-z|}(u_\varepsilon^{(2)}(z))^2\,dz=O(
\varepsilon+\varepsilon^{-1}||\varphi_\varepsilon^{(2)}||^2_\varepsilon),
\end{equation}
\begin{align}
{}&\int_{\mathbb{R}^2}\frac{1}{|x-y|}(u_\varepsilon^{(1)}(y)+u_\varepsilon^{(2)}(y))
\xi_\varepsilon(y)\,dy
=O\Big(\sum_{i=1}^2\int_{\mathbb{R}^2}\frac{1}{|x-y|}|u_\varepsilon^{(i)}(y)\xi_\varepsilon(y)|\,dy\Big)\nonumber\\
={}&O\Big(\sum_{i=1}^2\sum_{j=1}^k\int_{\mathbb{R}^2}\frac{1}{|x-y|}|U^j_{\varepsilon,y_\varepsilon^{j(i)}}(y)\xi_\varepsilon(y)|\,dy\Big)
+O\Big(\sum_{i=1}^2\int_{\mathbb{R}^2}\frac{1}{|x-y|}|\varphi_\varepsilon^{(i)}(y)\xi_\varepsilon(y)|\,dy\Big)\nonumber\\
={}&O\Big(\sum_{i=1}^2\sum_{j=1}^k(\int_{\mathbb{R}^2}\frac{|U^j_{\varepsilon,y_\varepsilon^{j(i)}}(y)|^2}{|x-y|^2}\,dy)
^{\frac{1}{2}}(\int_{\mathbb{R}^2}|\xi_\varepsilon(y)|^2\,dy)^{\frac{1}{2}}\Big)\label{151}\\
{}&+O\Big(\sum_{i=1}^2(\int_{\mathbb{R}^2}\frac{|\varphi_\varepsilon^{(i)}(y)|^2}{|x-y|}\,dy)^{\frac{1}{2}}
(\int_{\mathbb{R}^2}\frac{|\xi_\varepsilon(y)|^2}{|x-y|}\,dy)^{\frac{1}{2}}\Big)\nonumber\\
={}&O\Big(1+\sum_{i=1}^2\varepsilon^{-1}||\varphi_\varepsilon^{(i)}||_\varepsilon\Big)||\xi_\varepsilon||_\varepsilon,\nonumber
\end{align}
and
\begin{equation}\label{152}
\begin{aligned}
\int_{\mathbb{R}^2}u_\varepsilon^{(2)}(x)\xi_\varepsilon(x)\,dx
&=O\Big(\sum_{j=1}^k\int_{\mathbb{R}^2}|U^j_{\varepsilon,y_\varepsilon^{j(2)}}(x)\xi_\varepsilon(x)|\,dx\Big)
+O\Big(\int_{\mathbb{R}^2}|\varphi_\varepsilon^{(2)}(x)\xi_\varepsilon(x)|\,dx\Big)\\
&=O(\varepsilon+
||\varphi_\varepsilon^{(2)}||_\varepsilon)||\xi_\varepsilon||_\varepsilon.
\end{aligned}
\end{equation}
Then it follows from \eqref{109} that
\begin{align*}
{}&\int_{\mathbb{R}^2}E_{1}(x)\xi_\varepsilon\,dx\\
={}&O\Big(\int_{\mathbb{R}^2}u_\varepsilon^{(2)}(x)\xi_\varepsilon(x)
\Big(\int_{\mathbb{R}^2}\frac{1}{|x-y|}((u_\varepsilon^{(1)}(y)+u_\varepsilon^{(2)}(y))
\xi_\varepsilon(y))\Big(\int_{\mathbb{R}^2}\frac{1}{|y-z|}(u_\varepsilon^{(2)}(z))^2\,dz\Big)\,dy\Big)\,dx\Big)\\
={}&O\Big((\varepsilon+\varepsilon^{-1}||\varphi_\varepsilon^{(2)}||^2_\varepsilon)
\Big(1+\sum_{i=1}^2\varepsilon^{-1}||\varphi_\varepsilon^{(i)}||_\varepsilon\Big)||\xi_\varepsilon||_\varepsilon
(\varepsilon+
||\varphi_\varepsilon^{(2)}||_\varepsilon)||\xi_\varepsilon||_\varepsilon\Big)=o(1)||\xi_\varepsilon||^2_\varepsilon,
\end{align*}
\begin{align*}
{}&\int_{\mathbb{R}^2}E_{2}(x)\xi_\varepsilon\,dx\\
={}&O\Big(\int_{\mathbb{R}^2}u_\varepsilon^{(2)}(x)\xi_\varepsilon(x)\Big(\int_{\mathbb{R}^2}\frac{1}{|x-y|}(u_\varepsilon^{(1)}(y))^2
\Big(\int_{\mathbb{R}^2}\frac{1}{|y-z|}(u_\varepsilon^{(1)}(z)+u_\varepsilon^{(2)}(z))\xi_\varepsilon(z)\,dz\Big)\,dy\Big)\,dx\Big)\\
={}&O\Big(\Big(1+\sum_{i=1}^2\varepsilon^{-1}||\varphi_\varepsilon^{(i)}||_\varepsilon\Big)||\xi_\varepsilon||_\varepsilon
(\varepsilon+\varepsilon^{-1}||\varphi_\varepsilon^{(1)}||^2_\varepsilon)(\varepsilon+
||\varphi_\varepsilon^{(2)}||_\varepsilon)||\xi_\varepsilon||_\varepsilon\Big)=o(1)||\xi_\varepsilon||^2_\varepsilon,
\end{align*}
and
\begin{align*}
{}&\int_{\mathbb{R}^2}E_{3}(x)\xi_\varepsilon\,dx\\
={}&O\Big(\int_{\mathbb{R}^2}u_\varepsilon^{(2)}(x)\xi_\varepsilon(x)\Big(\sum_{i=1}^2\int_{\mathbb{R}^2}\frac{1}{|x-y|}((u_\varepsilon^{(i)}(y))^2\,dy\Big)
\Big(\int_{\mathbb{R}^2}\frac{1}{|x-z|}(u_\varepsilon^{(1)}(z)+u_\varepsilon^{(2)}(z))\xi_\varepsilon(z)\,dz\Big)\,dx\Big)\\
={}&O\Big(\Big(\sum_{i=1}^2(\varepsilon+\varepsilon^{-1}||\varphi_\varepsilon^{(i)}||^2_\varepsilon)\Big)\Big(1+\sum_{i=1}^2
\varepsilon^{-1}||\varphi_\varepsilon^{(i)}||_\varepsilon\Big)||\xi_\varepsilon||_\varepsilon(\varepsilon+
||\varphi_\varepsilon^{(2)}||_\varepsilon)||\xi_\varepsilon||_\varepsilon\Big)=o(1)||\xi_\varepsilon||^2_\varepsilon.
\end{align*}
The rest terms $\int_{\mathbb{R}^2}E_{i}(x)\xi_\varepsilon\,dx(i=4,5,6)$ can be estimated similarly. Hence
\begin{equation}\label{154}\nonumber
\int_{\mathbb{R}^2}E(x)\xi_\varepsilon\,dx=o(1)||\xi_\varepsilon||^2_\varepsilon.
\end{equation}
From the above analysis we obtain \eqref{146}.
\end{proof}

\begin{lem}\label{prop10}
Suppose $\varphi_\varepsilon(x)$ is derived as in Proposition \ref{prop2} with $||\varphi_\varepsilon||_\varepsilon=o(\varepsilon)$. Then there exist $C>0$ and $\sigma>0$ such that
\begin{equation}\label{178}
|\varphi_\varepsilon(x)|\leq C\sum_{j=1}^ke^{-\frac{\sigma|x-y^j_\varepsilon|}{\varepsilon}},\quad \forall x\in\mathbb{R}^2.
\end{equation}
\end{lem}

\noindent
\begin{proof}
Let $u^i_\varepsilon(x)=u_\varepsilon(\varepsilon x+y_\varepsilon^{i})$ and $\varphi^i_\varepsilon(x)=\varphi_\varepsilon(\varepsilon x+y_\varepsilon^{i})$, then $\varphi^i_\varepsilon$ satisfies
\begin{equation}\nonumber
-\Delta\varphi^i_\varepsilon+V(\varepsilon x+y_\varepsilon^{i})\varphi^i_\varepsilon=f_\varepsilon(x),
\end{equation}
where
\begin{align*}
f_\varepsilon(x):={}&
\Big(\sum_{j=1}^kU^j_{\varepsilon,y^j_\varepsilon}(\varepsilon x+y_\varepsilon^{i})+\varphi^i_\varepsilon\Big)^{p-1}-\Big(\sum_{j=1}^kU^j_{\varepsilon,y^j_\varepsilon}(\varepsilon x+y_\varepsilon^{i})\Big)^{p-1}\\
{}&+\Big(\sum_{j=1}^kU^j_{\varepsilon,y^j_\varepsilon}(\varepsilon x+y_\varepsilon^{i})\Big)^{p-1}-
\sum_{j=1}^k\Big(U^j_{\varepsilon,y^j_\varepsilon}(\varepsilon x+y_\varepsilon^{i})\Big)^{p-1}\\
{}&-\sum_{j=1}^k(V(\varepsilon x+y_\varepsilon^{i})-V(a^j))U^j_{\varepsilon,y^j_\varepsilon}(\varepsilon x+y_\varepsilon^{i})\\
{}&-\varepsilon^2(A_0(u^i_{\varepsilon})+A_1^2(u^i_{\varepsilon})+A_2^2(
u^i_{\varepsilon}))u^i_{\varepsilon}.
\end{align*}
Since $||\varphi_\varepsilon||_\varepsilon=o(\varepsilon)$, we have $||\varphi^i_\varepsilon||_{H^1(\mathbb{R}^2)}=o(1)$. Then by Moser iteration we obtain
\begin{equation}\nonumber
||\varphi_\varepsilon||_{L^\infty(\mathbb{R}^2)}=||\varphi^i_\varepsilon||_{L^\infty(\mathbb{R}^2)}=o(1),
\end{equation}
which leads to \eqref{178} by applying the comparison theorem.
\end{proof}

\begin{lem}\label{prop4}
For any fixed small $d>0$, there holds
\begin{equation}\label{155}
E(x)=o(1)\ \text{in}\ B_{d}(y_\varepsilon^{i(1)}).
\end{equation}
\end{lem}

\noindent
\begin{proof}
By Proposition \ref{prop3}, we deduce $|y_\varepsilon^{i(1)}-y_\varepsilon^{i(2)}|=o(\varepsilon)$.
Then for $x\in B_{d}(y_\varepsilon^{i(1)})$, it follows from \eqref{114} and Lemma \ref{prop10} that
\begin{equation}\label{156}
|u_\varepsilon^{(2)}(x)|\leq C\sum_{j=1}^ke^{-\frac{\sigma|x-y^{j(2)}_\varepsilon|}{\varepsilon}}\leq C
e^{-\frac{\sigma|x-y^{j(1)}_\varepsilon|}{\varepsilon}+o(1)}\leq C.
\end{equation}
Combining \eqref{109}, \eqref{146}, \eqref{150}, \eqref{151} and \eqref{156} yields
\begin{equation}\label{157}\nonumber
\begin{aligned}
\big|E_{1}(x)\big|
\leq{}&C|u_\varepsilon^{(2)}(x)|\Big|\int_{\mathbb{R}^2}\frac{1}{|x-y|}(u_\varepsilon^{(1)}(y)+u_\varepsilon^{(2)}(y))
\xi_\varepsilon(y)
\Big(\int_{\mathbb{R}^2}\frac{1}{|y-z|}(u_\varepsilon^{(2)}(z))^2\,dz\Big)\,dy\Big|\\
\leq{}&C(\varepsilon+\varepsilon^{-1}||\varphi_\varepsilon^{(2)}||^2_\varepsilon)
\Big(1+\sum_{i=1}^2\varepsilon^{-1}||\varphi_\varepsilon^{(i)}||_\varepsilon\Big)||\xi_\varepsilon||_\varepsilon=o(1),
\end{aligned}
\end{equation}
\begin{align*}
\big|E_{2}(x)\big|
\leq{}&C|u_\varepsilon^{(2)}(x)|\Big|\int_{\mathbb{R}^2}\frac{1}{|x-y|}(u_\varepsilon^{(1)}(y))^2
\Big(\int_{\mathbb{R}^2}\frac{1}{|y-z|}(u_\varepsilon^{(1)}(z)+u_\varepsilon^{(2)}(z))\xi_\varepsilon(z)\,dz\Big)\,dy\Big|\\
\leq{}&C
\Big(1+\sum_{i=1}^2\varepsilon^{-1}||\varphi_\varepsilon^{(i)}||_\varepsilon\Big)||\xi_\varepsilon||_\varepsilon
(\varepsilon+\varepsilon^{-1}||\varphi_\varepsilon^{(1)}||^2_\varepsilon)=o(1),
\end{align*}
and
\begin{align*}
\big|E_{3}(x)\big|
\leq{}&C|u_\varepsilon^{(2)}(x)|
\Big|\sum_{i=1}^2\int_{\mathbb{R}^2}\frac{1}{|x-y|}((u_\varepsilon^{(i)}(y))^2\,dy\cdot
\int_{\mathbb{R}^2}\frac{1}{|x-z|}(u_\varepsilon^{(1)}(z)+u_\varepsilon^{(2)}(z))\xi_\varepsilon(z)\,dz\Big|\\
\leq{}&C\sum_{i=1}^2(\varepsilon+\varepsilon^{-1}||\varphi_\varepsilon^{(i)}||^2_\varepsilon)
\Big(1+\sum_{i=1}^2\varepsilon^{-1}||\varphi_\varepsilon^{(i)}||_\varepsilon\Big)||\xi_\varepsilon||_\varepsilon=o(1).
\end{align*}
Analogously, $E_4(x),E_5(x)$ and $E_6(x)$ can be estimated. Hence we obtain \eqref{155}.
\end{proof}

\begin{lem}\label{prop11}
For large $R>0$, there holds
\begin{equation}
E(x)=o(1)\ \text{in}\ \mathbb{R}^2\backslash\cup_{i=1}^kB_{R\varepsilon}(y_\varepsilon^{i(1)}).
\end{equation}
\end{lem}

\noindent
\begin{proof}
Given any small $\gamma>0$, for $x\in\mathbb{R}^2\backslash\cup_{i=1}^kB_{R\varepsilon}(y_\varepsilon^{i(1)})$, by \eqref{114} and Lemma \ref{prop10}, we have
\begin{equation}
|u_\varepsilon^{(2)}(x)|\leq C\sum_{j=1}^ke^{-\frac{\sigma|x-y^{j(2)}_\varepsilon|}{\varepsilon}}\leq C\sum_{j=1}^k
e^{-\frac{\sigma|x-y^{j(1)}_\varepsilon|}{\varepsilon}+o(1)}\leq \gamma.
\end{equation}
Similar to the proof of Lemma \ref{prop4} we complete the proof.
\end{proof}

\begin{prop}\label{prop5}
Let $\xi_\varepsilon^i(x)=\xi_\varepsilon(\varepsilon x+y_\varepsilon^{i(1)})$ for $i=1,\cdots,k$. Then there exist $d_{j,i}\in\mathbb{R},\ i=1,\cdots,k,\ j=1,2,$ such that, by taking a subsequence necessarily, there holds
\begin{equation}\label{140}
\xi_\varepsilon^i(x)\to\sum_{j=1}^2d_{j,i}\frac{\partial U^i(x)}{\partial x_j}\quad\text{in}\ C^1_{loc}(\mathbb{R}^2),
\end{equation}
as $\varepsilon\to0$.
\end{prop}

\noindent
\begin{proof}
Calculating directly, we deduce
\begin{equation}\label{141}\nonumber
\begin{aligned}
{}&-\Delta\xi^i_\varepsilon+V(\varepsilon x+{y^{i(1)}_\varepsilon})\xi^i_\varepsilon+\varepsilon^2(A_0(u_\varepsilon^{(1)})+ A_1^2(u_\varepsilon^{(1)})+A_2^2(u_\varepsilon^{(1)}))\xi^i_\varepsilon\\
={}&C_\varepsilon(\varepsilon x+{y^{i(1)}_\varepsilon})\xi^i_\varepsilon+E(\varepsilon x+{y^{i(1)}_\varepsilon}).
\end{aligned}
\end{equation}
Based on the fact $||\xi_\varepsilon^i||_{L^\infty(\mathbb{R}^2)}=1$, by the elliptic regularity theory, we have
$\xi_\varepsilon^i\in C^{1,\alpha}_{loc}(\mathbb{R}^2)$ and $||\xi_\varepsilon^i||_{C^{1,\alpha}_{loc}(\mathbb{R}^2)}\leq C$ for some $\alpha\in(0,1)$. Then up to a subsequence, we can assume $\xi_\varepsilon^i\to\xi^i$ in $C^{1}_{loc}(\mathbb{R}^2)$. For any given $\psi\in C^\infty_0(\mathbb{R}^2)$, we have
\begin{equation}\label{142}\nonumber
\begin{aligned}
&\quad\int_{\mathbb{R}^2}(\nabla\xi^i_\varepsilon\nabla\psi+V(\varepsilon x+{y^{i(1)}_\varepsilon})\xi^i_\varepsilon\psi)\,dx+\varepsilon^2\int_{\mathbb{R}^2}
(A_0(u_\varepsilon^{(1)})+ A_1^2(u_\varepsilon^{(1)})+A_2^2(u_\varepsilon^{(1)}))\xi^i_\varepsilon\psi\,dx\\
&=\int_{\mathbb{R}^2}(C_\varepsilon(\varepsilon x+{y^{i(1)}_\varepsilon})\xi^i_\varepsilon\psi+E(\varepsilon x+{y^{i(1)}_\varepsilon})\psi)\,dx.
\end{aligned}
\end{equation}
Combining \eqref{120}, \eqref{121} and Lemma \ref{prop4} yields
\begin{equation}\label{143}
\int_{\mathbb{R}^2}(\nabla\xi^i_\varepsilon\nabla\psi+V(\varepsilon x+{y^{i(1)}_\varepsilon})\xi^i_\varepsilon\psi)\,dx=\int_{\mathbb{R}^2}C_\varepsilon(\varepsilon x+{y^{i(1)}_\varepsilon})\xi^i_\varepsilon\psi\,dx+o(1)||\psi||_{H^1(\mathbb{R}^2)}.
\end{equation}
Similar to Cao et al. \cite{ref4}, it turns out that
\begin{equation}\label{144}\nonumber
\int_{\mathbb{R}^2}C_\varepsilon(\varepsilon x+{y^{i(1)}_\varepsilon})\xi^i_\varepsilon\psi\,dx\to
\int_{\mathbb{R}^2}(p-1)(U^i)^{p-2}\xi^i\psi\,dx.
\end{equation}
Letting $\varepsilon\to0$ in \eqref{143}, we find $\xi^i$ satisfies
\begin{equation}\label{145}
-\Delta\xi^i+V(a^i)\xi^i=(p-1)(U^i)^{p-2}\xi^i.
\end{equation}
Then the nondegeneracy of \eqref{145} implies \eqref{140}.
\end{proof}

\begin{prop}\label{prop6}
Let $d_{j,i}$ be defined as in Proposition \ref{prop5}. Then
\begin{equation}\label{158}
d_{j,i}=0,\quad j=1,2.
\end{equation}
\end{prop}

\noindent
\begin{proof}
Applying the Pohozaev identity \eqref{P} to $u_\varepsilon^{(1)}$ and $u_\varepsilon^{(2)}$ with $\Omega=B_d(y_\varepsilon^{i(1)})$ for some small constant $d>0$, we have
\begin{equation}\label{159}
\begin{aligned}
{}&\int_{B_d(y_\varepsilon^{i(1)})}\frac{\partial V}{\partial x_k}(u_\varepsilon^{(1)}+u_\varepsilon^{(2)})\xi_\varepsilon\,dx\\
={}&\int_{\partial B_d(y_\varepsilon^{i(1)})}\Big(\varepsilon^2\nabla(u_\varepsilon^{(1)}+u_\varepsilon^{(2)})\nabla\xi_\varepsilon
+V(x)(u_\varepsilon^{(1)}+u_\varepsilon^{(2)})\xi_\varepsilon\Big)\nu_k\,d\sigma\\
{}&-2\int_{\partial B_d(y_\varepsilon^{i(1)})}\Big(\varepsilon^2\Big(\frac{\partial u_\varepsilon^{(1)}}{\partial\nu}
\frac{\partial \xi_\varepsilon}{\partial x_k}+\frac{\partial \xi_\varepsilon}{\partial\nu}\frac{\partial u_\varepsilon^{(2)}}{\partial x_k}\Big)
+F_1(x)\xi_\varepsilon\nu_k\Big)\,d\sigma\\
{}&-2\int_{\partial B_d(y_\varepsilon^{i(1)})}F_2(x)\,d\sigma+\int_{\partial B_d(y_\varepsilon^{i(1)})}F_3(x)\nu_k\,d\sigma,
\end{aligned}
\end{equation}
where
\begin{equation}\label{160}\nonumber
F_1(x)=\int_0^1(tu_\varepsilon^{(1)}+(1-t)u_\varepsilon^{(2)})^{p-1}\,dt,
\end{equation}
\begin{align*}\label{161}
F_2(x)={}&(A_1(u_\varepsilon^{(1)}(x))\nu_1(x)+
A_2(u_\varepsilon^{(1)}(x))\nu_2(x))A_k(u_\varepsilon^{(1)}(x))(u_\varepsilon^{(1)}(x)+u_\varepsilon^{(2)}(x)
)\xi_\varepsilon(x)\\
{}&+(A_1(u_\varepsilon^{(2)}(x))\nu_1(x)+
A_2(u_\varepsilon^{(2)}(x))\nu_2(x))(u_\varepsilon^{(2)}(x))^2\\
{}&\quad\cdot\Big(\frac{(-1)^k}{4\pi}\int_{\mathbb{R}^2}
\frac{x_{k+(-1)^{k+1}}-y_{k+(-1)^{k+1}}}{|x-y|^2}(u_\varepsilon^{(1)}(y)+u_\varepsilon^{(2)}(y))\xi_\varepsilon(y)\,dy\Big)\\
{}&+A_k(u_\varepsilon^{(1)}(x))(u_\varepsilon^{(2)}(x))^2\nu_1(x)\Big(
-\frac{1}{4\pi}\int_{\mathbb{R}^2}\frac{x_2-y_2}{|x-y|^2}(u_\varepsilon^{(1)}(y)+u_\varepsilon^{(2)}(y))
\xi_\varepsilon(y)\,dy\Big)\\
{}&+A_k(u_\varepsilon^{(1)}(x))(u_\varepsilon^{(2)}(x))^2\nu_2(x)\Big(
\frac{1}{4\pi}\int_{\mathbb{R}^2}\frac{x_1-y_1}{|x-y|^2}(u_\varepsilon^{(1)}(y)+u_\varepsilon^{(2)}(y))
\xi_\varepsilon(y)\,dy\Big),
\end{align*}
and
\begin{align*}
F_3(x)={}&(A_0(u_\varepsilon^{(1)}(x))+ A_1^2(u_\varepsilon^{(1)}(x))+A_2^2(u_\varepsilon^{(1)}(x))
)(u_\varepsilon^{(1)}(x)+u_\varepsilon^{(2)}(x))\xi_\varepsilon(x)\\
{}&-\frac{(u_\varepsilon^{(2)}(x))^2}{8\pi^2}\int_{\mathbb{R}^2}\frac{x_1-y_1}{|x-y|^2}
(u_\varepsilon^{(1)}(y))^2\Big(\int_{\mathbb{R}^2}\frac{y_1-z_1}{|y-z|^2}(u_\varepsilon^{(1)}(z)
+u_\varepsilon^{(2)}(z))\xi_\varepsilon(z)\,dz\Big)\,dy\\
{}&-\frac{(u_\varepsilon^{(2)}(x))^2}{8\pi^2}\int_{\mathbb{R}^2}\frac{x_1-y_1}{|x-y|^2}
(u_\varepsilon^{(1)}(y)+u_\varepsilon^{(2)}(y))\xi_\varepsilon(y)\Big(\int_{\mathbb{R}^2}
\frac{y_1-z_1}{|y-z|^2}(u_\varepsilon^{(2)}(z))^2\,dz\Big)\,dy\\
{}&-\frac{(u_\varepsilon^{(2)}(x))^2}{8\pi^2}\int_{\mathbb{R}^2}\frac{x_2-y_2}{|x-y|^2}
(u_\varepsilon^{(1)}(y))^2\Big(\int_{\mathbb{R}^2}\frac{y_2-z_2}{|y-z|^2}(u_\varepsilon^{(1)}(z)
+u_\varepsilon^{(2)}(z))\xi_\varepsilon(z)\,dz\Big)\,dy\\
{}&-\frac{(u_\varepsilon^{(2)}(x))^2}{8\pi^2}\int_{\mathbb{R}^2}\frac{x_2-y_2}{|x-y|^2}
(u_\varepsilon^{(1)}(y)+u_\varepsilon^{(2)}(y))\xi_\varepsilon(y)\Big(\int_{\mathbb{R}^2}
\frac{y_2-z_2}{|y-z|^2}(u_\varepsilon^{(2)}(z))^2\,dz\Big)\,dy\\
{}&+(u_\varepsilon^{(2)}(x))^2(A_1(u_\varepsilon^{(1)}(x))+A_1(u_\varepsilon^{(2)}(x)))\Big(-\frac{1}{4\pi}\int_{\mathbb{R}^2}
\frac{x_2-y_2}{|x-y|^2}(u_\varepsilon^{(1)}(y)+u_\varepsilon^{(2)}(y))\xi_\varepsilon(y)\,dy\Big)\\
{}&+(u_\varepsilon^{(2)}(x))^2(A_2(u_\varepsilon^{(1)}(x))+A_2(u_\varepsilon^{(2)}(x)))\Big(\frac{1}{4\pi}\int_{\mathbb{R}^2}
\frac{x_1-y_1}{|x-y|^2}(u_\varepsilon^{(1)}(y)+u_\varepsilon^{(2)}(y))\xi_\varepsilon(y)\,dy\Big).
\end{align*}
By using \eqref{109}, \eqref{116}-\eqref{119} and \eqref{146} and choosing $\gamma>0$ sufficiently large, we derive
\begin{equation}\label{163}
\begin{aligned}
{}&\int_{\partial B_d(y_\varepsilon^{i(1)})}\Big(\varepsilon^2\nabla(u_\varepsilon^{(1)}+u_\varepsilon^{(2)})\nabla\xi_\varepsilon
+V(x)(u_\varepsilon^{(1)}+u_\varepsilon^{(2)})\xi_\varepsilon\Big)\nu_k\,d\sigma\\
={}&O\Big(\sum_{l=1}^2(\varepsilon^\gamma+||\varphi_\varepsilon^{(l)}||^2_\varepsilon)^\frac{1}{2}||\xi_\varepsilon||_\varepsilon\Big)
=O(\varepsilon^{m+2}),
\end{aligned}
\end{equation}
\begin{equation}\label{164}
\begin{aligned}
{}&\int_{\partial B_d(y_\varepsilon^{i(1)})}\varepsilon^2\Big(\frac{\partial u_\varepsilon^{(1)}}{\partial\nu}
\frac{\partial \xi_\varepsilon}{\partial x_k}+\frac{\partial \xi_\varepsilon}{\partial\nu}\frac{\partial u_\varepsilon^{(2)}}{\partial x_k}\Big)\,d\sigma\\
={}&O\Big(\sum_{l=1}^2(\varepsilon^\gamma+||\varphi_\varepsilon^{(l)}||^2_\varepsilon)^\frac{1}{2}
||\xi_\varepsilon||_\varepsilon\Big)=O(\varepsilon^{m+2}),
\end{aligned}
\end{equation}
and
\begin{equation}\label{165}
\begin{aligned}
{}&\int_{\partial B_d(y_\varepsilon^{i(1)})}F_1(x)\xi_\varepsilon\nu_k\,d\sigma\\
={}&O\Big(\sum_{l=1}^2\int_{\partial B_d(y_\varepsilon^{i(1)})}|u_\varepsilon^{(l)}|^{p-1}\xi_\varepsilon\,d\sigma\Big)\\
={}&O\Big(\sum_{l=1}^2\Big(\int_{\partial B_d(y_\varepsilon^{i(1)})}|u_\varepsilon^{(l)}|^{p}\,d\sigma\Big)^{\frac{p-1}{p}}
||\xi_\varepsilon||_{L^p(\mathbb{R}^2)}\Big)\\
={}&O\Big(\sum_{l=1}^2(\varepsilon^\gamma+||\varphi_\varepsilon^{(l)}||^2_\varepsilon)^\frac{p-1}{p}
(\varepsilon^{2-p}||\xi_\varepsilon||^p_\varepsilon)^{\frac{1}{p}}\Big)\\
={}&O(\varepsilon^{2m+2-\frac{2m}{p}}).
\end{aligned}
\end{equation}
Also, by using \eqref{109}, \eqref{118}, \eqref{120}, \eqref{146} and \eqref{151}, we obtain
\begin{equation}\label{166}\nonumber
\begin{aligned}
{}&\int_{\partial B_d(y_\varepsilon^{i(1)})}(A_1(u_\varepsilon^{(1)}(x))\nu_1(x)+
A_2(u_\varepsilon^{(1)}(x))\nu_2(x))A_k(u_\varepsilon^{(1)}(x))(u_\varepsilon^{(1)}(x)+u_\varepsilon^{(2)}(x)
)\xi_\varepsilon(x)\,d\sigma\\
={}&O\Big((\varepsilon^2+\varepsilon^{-2}||\varphi_\varepsilon^{(1)}||^4_\varepsilon)\sum_{l=1}^2\int_{\partial B_d(y_\varepsilon^{i(1)})}u_\varepsilon^{(l)}(x)\xi_\varepsilon(x)\,d\sigma\Big)\\
={}&O\Big((\varepsilon^2+\varepsilon^{-2}||\varphi_\varepsilon^{(1)}||^4_\varepsilon)\sum_{l=1}^2
(\varepsilon^\gamma+||\varphi_\varepsilon^{(l)}||^2_\varepsilon)^\frac{1}{2}||\xi_\varepsilon||_{\varepsilon}\Big)
=O(\varepsilon^{m+4}),
\end{aligned}
\end{equation}
and
\begin{align*}
{}&\int_{\partial B_d(y_\varepsilon^{i(1)})}(A_1(u_\varepsilon^{(2)}(x))\nu_1(x)+
A_2(u_\varepsilon^{(2)}(x))\nu_2(x))(u_\varepsilon^{(2)}(x))^2\\
{}&\quad\cdot\Big(\frac{(-1)^k}{4\pi}\int_{\mathbb{R}^2}
\frac{x_{k+(-1)^{k+1}}-y_{k+(-1)^{k+1}}}{|x-y|^2}(u_\varepsilon^
{(1)}(y)+u_\varepsilon^{(2)}(y))\xi_\varepsilon(y)\,dy\Big)\,d\sigma\\
={}&O\Big((\varepsilon+\varepsilon^{-1}||\varphi_\varepsilon^{(2)}||^2_\varepsilon)
\Big(1+\sum_{l=1}^2\varepsilon^{-1}||\varphi_\varepsilon^{(l)}||_\varepsilon\Big)||\xi_\varepsilon||_\varepsilon\int_{\partial B_d(y_\varepsilon^{i(1)})}(u_\varepsilon^{(2)}(x))^2\,d\sigma\Big)\\
={}&O\Big((\varepsilon+\varepsilon^{-1}||\varphi_\varepsilon^{(2)}||^2_\varepsilon)
\Big(1+\sum_{l=1}^2\varepsilon^{-1}||\varphi_\varepsilon^{(l)}||_\varepsilon\Big)||\xi_\varepsilon||_\varepsilon
(\varepsilon^\gamma+||\varphi_\varepsilon^{(2)}||^2_\varepsilon)\Big)=O(\varepsilon^{2m+4}).
\end{align*}

Similarly, we can verify that
\begin{equation}\nonumber
\int_{\partial B_d(y_\varepsilon^{i(1)})}A_k(u_\varepsilon^{(1)}(x))(u_\varepsilon^{(2)}(x))^2\nu_1(x)\Big(
-\frac{1}{4\pi}\int_{\mathbb{R}^2}\frac{x_2-y_2}{|x-y|^2}(u_\varepsilon^{(1)}(y)+u_\varepsilon^{(2)}(y))
\xi_\varepsilon(y)\,dy\Big)\,d\sigma=O(\varepsilon^{2m+4}),
\end{equation}
\begin{equation}\nonumber
\int_{\partial B_d(y_\varepsilon^{i(1)})}A_k(u_\varepsilon^{(1)}(x))(u_\varepsilon^{(2)}(x))^2\nu_2(x)\Big(
\frac{1}{4\pi}\int_{\mathbb{R}^2}\frac{x_1-y_1}{|x-y|^2}(u_\varepsilon^{(1)}(y)+u_\varepsilon^{(2)}(y))
\xi_\varepsilon(y)\,dy\Big)\,d\sigma=O(\varepsilon^{2m+4}).
\end{equation}
Hence
\begin{equation}\label{168}
\int_{\partial B_d(y_\varepsilon^{i(1)})}F_2(x)\,d\sigma=O(\varepsilon^{m+4}).
\end{equation}
By the same token, we have
\begin{equation}\label{169}
\int_{\partial B_d(y_\varepsilon^{i(1)})}F_3(x)\nu_k\,d\sigma=O(\varepsilon^{m+4}).
\end{equation}
Combining \eqref{163}-\eqref{169} yields
\begin{equation}\label{170}
\text{RHS of \eqref{159}}=O(\varepsilon^{m+2}).
\end{equation}

On the other hand, it follows from $(V_3)$ that
\begin{equation}\label{171}
\begin{aligned}
\text{LHS of \eqref{159}}={}&\sum_{l=1}^2\int_{B_d(y_\varepsilon^{i(1)})}mb_{k,i}|x_k-a^i_k|^{m-2}
(x_k-a^i_k)u_\varepsilon^{(l)}\xi_\varepsilon\,dx\\
{}&+O\Big(\sum_{l=1}^2\int_{B_d(y_\varepsilon^{i(1)})}|x-a^i|^{m}u_\varepsilon^{(l)}\xi_\varepsilon\,dx\Big).
\end{aligned}
\end{equation}
Then using the fact that $|\xi_\varepsilon|\leq1$ and Proposition \ref{prop3}, we can deduce that
\begin{align}
{}&\int_{B_d(y_\varepsilon^{i(1)})}mb_{k,i}|x_k-a^i_k|^{m-2}
(x_k-a^i_k)u_\varepsilon^{(l)}\xi_\varepsilon\,dx\nonumber\\
={}&\int_{B_d(y_\varepsilon^{i(1)})}mb_{k,i}|x_k-a^i_k|^{m-2}
(x_k-a^i_k)U^i_{\varepsilon,y_\varepsilon^{i(l)}}\xi_\varepsilon\,dx\nonumber\\
{}&+O\Big(\sum_{j\neq i}^k\int_{B_d(y_\varepsilon^{i(1)})}|x-a^i|^{m-1}U^j_{\varepsilon,y_\varepsilon^{j(l)}}\xi_\varepsilon\,dx
+\int_{B_d(y_\varepsilon^{i(1)})}|x-a^i|^{m-1}\varphi_\varepsilon^{(l)}\xi_\varepsilon\,dx\Big)\nonumber\\
={}&\varepsilon^2\int_{B_{\frac{d}{\varepsilon}}(0)}mb_{k,i}|\varepsilon x_k+y_{\varepsilon,k}^{i(1)}-a^i_k|^{m-2}
(\varepsilon x_k+y_{\varepsilon,k}^{i(1)}-a^i_k)U^i(x+\frac{y_{\varepsilon}^{i(1)}-y_{\varepsilon}^{i(l)}}{\varepsilon})\xi^i_\varepsilon(x)\,dx\label{172}\\
{}&+O(e^{-\frac{\tau}{\varepsilon}}+||\varphi_\varepsilon^{(l)}||_\varepsilon||\xi_\varepsilon||_\varepsilon)\nonumber\\
={}&\varepsilon^2\int_{B_{\frac{d}{\varepsilon}}(0)}mb_{k,i}|\varepsilon x_k+y_{\varepsilon,k}^{i(1)}-a^i_k|^{m-2}
(\varepsilon x_k+y_{\varepsilon,k}^{i(1)}-a^i_k)U^i(x+\frac{y_{\varepsilon}^{i(1)}-y_{\varepsilon}^{i(l)}}{\varepsilon})\xi^i_\varepsilon(x)\,dx\nonumber\\
{}&+O(\varepsilon^{m+2})\nonumber
\end{align}
and
\begin{equation}\label{173}
\begin{aligned}
{}&O\Big(\int_{B_d(y_\varepsilon^{i(1)})}|x-a^i|^{m}u_\varepsilon^{(l)}\xi_\varepsilon\,dx\Big)\\
={}&O\Big(\int_{B_d(y_\varepsilon^{i(1)})}|x-a^i|^{m}U^i_{\varepsilon,y_\varepsilon^{i(l)}}\xi_\varepsilon\,dx
+\int_{B_d(y_\varepsilon^{i(1)})}|x-a^i|^{m}\varphi_\varepsilon^{(l)}\xi_\varepsilon\,dx+e^{-\frac{\tau}{\varepsilon}}\Big)\\
={}&O\Big(\varepsilon^2\int_{B_{\frac{d}{\varepsilon}}(0)}|\varepsilon x+y_\varepsilon^{i(l)}-a^i|^mU^i(x+\frac{y_{\varepsilon}^{i(1)}-y_{\varepsilon}^{i(l)}}{\varepsilon})\,dx+
||\varphi_\varepsilon^{(l)}||_\varepsilon||\xi_\varepsilon||_\varepsilon
+e^{-\frac{\tau}{\varepsilon}}\Big)\\
={}&O(\varepsilon^{m+2}).
\end{aligned}
\end{equation}
Combining \eqref{170}-\eqref{173} yields
\begin{equation}\label{174}\nonumber
\int_{B_{\frac{d}{\varepsilon}}(0)}\Big|x_k+\frac{y_{\varepsilon,k}^{i(1)}-a^i_k}{\varepsilon}\Big|^{m-2}
\Big(x_k+\frac{y_{\varepsilon,k}^{i(1)}-a^i_k}{\varepsilon}\Big)U^i(x+\frac{y_{\varepsilon}^{i(1)}-y_{\varepsilon}^{i(l)}}{\varepsilon})
\xi^i_\varepsilon(x)\,dx=O(\varepsilon),
\end{equation}
which implies that
\begin{equation}\label{175}\nonumber
\sum_{j=1}^2d_{j,i}\int_{\mathbb{R}^2}|x_k|^{m-2}x_kU^i(x)\frac{\partial U^i}{\partial x_j}(x)\,dx=0.
\end{equation}
Hence $d_{j,i}=0$ for $j=1,2$.
\end{proof}

\begin{prop}\label{prop8}
For any fixed $R>0$, there holds
\begin{equation}\label{176}
||\xi_\varepsilon||_{L^\infty(\cup_{i=1}^kB_{R\varepsilon}(y_\varepsilon^{i(1)}))}=o(1).
\end{equation}
\end{prop}

\noindent
\begin{proof}
From Propositions \ref{prop5} and \ref{prop6}, we have for any fixed $R>0$,
\begin{equation}\nonumber
||\xi_\varepsilon^i||_{L^\infty(B_R(0))}=o(1),\quad i=1,\cdots,k,
\end{equation}
which implies that \eqref{176}.
\end{proof}

\begin{prop}\label{prop9}
For large $R>0$, there holds
\begin{equation}\label{177}
||\xi_\varepsilon||_{L^\infty(\mathbb{R}^2\backslash\cup_{i=1}^kB_{R\varepsilon}(y_\varepsilon^{i(1)}))}=o(1).
\end{equation}
\end{prop}

\noindent
\begin{proof}
Given any small $\gamma>0$, for $x\in\mathbb{R}^2\backslash\cup_{i=1}^kB_{R\varepsilon}(y_\varepsilon^{i(1)})$, it follows from \eqref{120}, \eqref{121}, Proposition \ref{prop3} and Lemma \ref{prop10} that
\begin{equation}
(A_0(u_\varepsilon^{(1)})+ A_1^2(u_\varepsilon^{(1)})+A_2^2(u_\varepsilon^{(1)}))\leq\gamma
\end{equation}
and
\begin{equation}
\begin{aligned}
|C_\varepsilon(x)|&\leq C(u_\varepsilon^{(1)}(x)+u_\varepsilon^{(2)}(x))^{p-2}\\
&\leq C\Big(\sum_{j=1}^k(e^{-\frac{\sigma|x-y^{j(1)}_\varepsilon|}{\varepsilon}}+
e^{-\frac{\sigma|x-y^{j(2)}_\varepsilon|}{\varepsilon}})\Big)^{p-2}\\
&\leq C\Big(\sum_{j=1}^k
e^{-\frac{\sigma|x-y^{j(1)}_\varepsilon|}{\varepsilon}+o(1)}\Big)^{p-2}\leq\gamma.
\end{aligned}
\end{equation}
So it turns out that
\begin{equation}
V(x)+A_0(u_\varepsilon^{(1)})+ A_1^2(u_\varepsilon^{(1)})+A_2^2(u_\varepsilon^{(1)})-C_\varepsilon(x)\geq\min_{x\in\mathbb{R}^2}V(x)-2\gamma>0.
\end{equation}
Therefore, by using Lemma \ref{prop11}, the equation of $\xi_\varepsilon$ and the maximum principle, we can get \eqref{177}.
\end{proof}

Finally, we give the proof of Theorem \ref{thm2}.

\noindent
\begin{proof}[\textbf{Proof of Theorem \ref{thm2}.}]
We argue by means of contraction. Suppose that $u_\varepsilon^{(1)}\neq u_\varepsilon^{(2)}$. From Propositions \ref{prop8} and \ref{prop9}, we have
$||\xi_\varepsilon||_{L^\infty(\mathbb{R}^2)}=o(1)$,
which is in contradiction with $||\xi_\varepsilon||_{L^\infty(\mathbb{R}^2)}=1$.
\end{proof}

\section*{Appendix}

\appendix
\renewcommand{\theequation}{A.\arabic{equation}}

\section{ The estimate of the energy functional}\label{seca}
In this section, we mainly estimate $I_{\varepsilon}(W_{\varepsilon,Y}).$
\begin{prop}\label{est-I}
Suppose that $p>2$ and $V(x)$ satisfies $(V_1)$ and $(V_2)$. Then for $\varepsilon>0$ sufficiently small and $Y\in D_{\delta}$, there holds that
\begin{equation}\label{eq2}
\begin{aligned}
I_{\varepsilon}(W_{\varepsilon,Y})=&\,
\Big(\frac{1}{2}-\frac{1}{p}\Big)\varepsilon^2\sum_{i=1}^k\int_{\mathbb{R}^2}(U^i(x))^p\,dx
+\frac{1}{2}\varepsilon^2\sum_{i=1}^k\int_{\mathbb{R}^2}(V (y^i)-V(a^i))(U^i(x))^2\,dx+O(\varepsilon^{2+\theta}).\\
\end{aligned}
\end{equation}
\end{prop}

\noindent
\begin{proof}
Since the integral $\int_{\mathbb{R}^2}\frac{1}{|x-y|}(U^i(y))^2\,dy$ is uniformly bounded for every $x\in\mathbb{R}^2$, we have
\begin{equation}\label{xiu1}
\int_{\mathbb{R}^2}\frac{1}{|x-y|}(U_{\varepsilon,y^i}^i(y))^2\,dy=
\varepsilon\int_{\mathbb{R}^2}\frac{1}{|z-\frac{x-y^i}{\varepsilon}|}(U^i(z))^2\,dz\leq C\varepsilon.
\end{equation}
This implies
\begin{equation}\label{eq7}
\begin{aligned}
{}&\Big|\frac{1}{2}\int_{\mathbb{R}^2}\Big(-\frac{1}{4\pi}\int_{\mathbb{R}^2}\frac{x_2-y_2}{|x-y|^2}
W_{\varepsilon,Y}^2(y)\,dy\Big)^2W_{\varepsilon,Y}^2(x)\,dx\Big|\\
\leq{}& C\int_{\mathbb{R}^2}\Big(\sum_{i=1}^k\int_{\mathbb{R}^2}\frac{1}{|x-y|}(U_{\varepsilon,y^i}^i(y))^2\,dy\Big)^2\Big(\sum_{i=1}^k U_{\varepsilon,y^i}^i(x)\Big)^2\,dx\\
\leq{}& C\varepsilon^2\int_{\mathbb{R}^2}\Big(\sum_{i=1}^k U_{\varepsilon,y^i}^i(x)\Big)^2\,dx\\
\leq{}& C\varepsilon^4.
\end{aligned}
\end{equation}
Analogously, we can estimate
\begin{equation}\label{eq8}
\Big|\frac{1}{2}\int_{\mathbb{R}^2}\Big(\frac{1}{4\pi}\int_{\mathbb{R}^2}\frac{x_1-y_1}{|x-y|^2}
W_{\varepsilon,Y}^2(y)\,dy\Big)^2W_{\varepsilon,Y}^2(x)\,dx\Big|\leq C\varepsilon^4.
\end{equation}
Noting that
\begin{equation}\nonumber
\int_{\mathbb{R}^2}(\varepsilon^2\nabla U_{\varepsilon,y^i}^i\nabla U_{\varepsilon,y^j}^j+V(a^i)U_{\varepsilon,y^i}^i U_{\varepsilon,y^j}^j)\,dx=\int_{\mathbb{R}^2}(U_{\varepsilon,y^i}^i)^{p-1}U_{\varepsilon,y^j}^j\,dx
\end{equation}
for each $i,j=1,\cdots,k$, we have
\begin{equation}\label{eq3}\nonumber
\begin{aligned}
{}&\frac{1}{2}\int_{\mathbb{R}^{2}}(\varepsilon^{2}|\nabla W_{\varepsilon,Y}|^{2}+V(x)|W_{\varepsilon,Y}|^2)\,dx-\frac{1}{p}\int_{\mathbb{R}^{2}}|W_{\varepsilon,Y}|^p\,dx\\
={}&\frac{1}{2}\sum_{i,j=1}^k\int_{\mathbb{R}^2}(\varepsilon^{2}\nabla U_{\varepsilon,y^i}^i\nabla U_{\varepsilon,y^j}^j+V(x)U_{\varepsilon,y^i}^i U_{\varepsilon,y^j}^j)\,dx-\frac{1}{p}\int_{\mathbb{R}^{2}}\Big(\sum_{i=1}^k U_{\varepsilon,y^i}^i\Big)^p\,dx\\
={}&\frac{1}{2}\sum_{i,j=1}^k\int_{\mathbb{R}^2}((V(x)-V(a^i))U_{\varepsilon,y^i}^i U_{\varepsilon,y^j}^j\,dx\\
{}&+\Big(\frac{1}{2}\sum_{i=1}^k\int_{\mathbb{R}^2}(U_{\varepsilon,y^i}^i)^p\,dx+\frac{1}{2}\sum_{i\neq j}\int_{\mathbb{R}^2}(U_{\varepsilon,y^i}^i)^{p-1} U_{\varepsilon,y^j}^j\,dx-\frac{1}{p}\int_{\mathbb{R}^{2}}\Big(\sum_{i=1}^k U_{\varepsilon,y^i}^i\Big)^p\,dx\Big)\\
=:{}&I_1+I_2.
\end{aligned}
\end{equation}
To estimate $I_1$, note that
\begin{align*}
I_1={}&\frac{1}{2}\sum_{i=1}^k\int_{\mathbb{R}^2}((V(x)-V(y^i))(U_{\varepsilon,y^i}^i)^2\,dx+
\frac{1}{2}\sum_{i=1}^k\int_{\mathbb{R}^2}((V(y^i)-V(a^i))(U_{\varepsilon,y^i}^i)^2\,dx\\
{}&+\frac{1}{2}\sum_{i\neq j}\int_{\mathbb{R}^2}((V(x)-V(y^i))U_{\varepsilon,y^i}^i U_{\varepsilon,y^j}^j\,dx+\frac{1}{2}\sum_{i\neq j}\int_{\mathbb{R}^2}((V(y^i)-V(a^i))U_{\varepsilon,y^i}^i U_{\varepsilon,y^j}^j\,dx.
\end{align*}
Since $U^i$ decays exponentially at infinity and $V(x)$ satisfies $(V_1)$ and $(V_2)$, we have
\begin{align*}
{}&\int_{\mathbb{R}^2}((V(x)-V(y^i))(U_{\varepsilon,y^i}^i)^2\,dx\\
={}&\int_{B_\delta(y^i)}((V(x)-V(y^i))(U_{\varepsilon,y^i}^i)^2\,dx+\int_{\mathbb{R}^2\backslash B_\delta(y^i)}((V(x)-V(y^i))(U_{\varepsilon,y^i}^i)^2\,dx\\
\leq{}&C\varepsilon^2\int_{B_{\frac{\delta}{\varepsilon}}(0)}|\epsilon x|^\theta(U^i(x))^2\,dx
+C\varepsilon^2\int_{\mathbb{R}^2\backslash B_{\frac{\delta}{\varepsilon}}(0)}(U^i(x))^2\,dx\leq C\varepsilon^{2+\theta},
\end{align*}
\begin{align*}
{}&\int_{\mathbb{R}^2}((V(x)-V(y^i))U_{\varepsilon,y^i}^i U_{\varepsilon,y^j}^j\,dx\\
={}&\int_{B_\delta(y^i)}((V(x)-V(y^i))U_{\varepsilon,y^i}^i U_{\varepsilon,y^j}^j\,dx+\int_{\mathbb{R}^2\backslash B_\delta(y^i)}((V(x)-V(y^i))U_{\varepsilon,y^i}^i U_{\varepsilon,y^j}^j\,dx\\
\leq{}&C\varepsilon^2\int_{B_{\frac{\delta}{\varepsilon}}(0)}|\epsilon x|^\theta U^i(x)U^j(x+\frac{y^i-y^j}{\varepsilon})\,dx+C\varepsilon^2\int_{\mathbb{R}^2\backslash B_{\frac{\delta}{\varepsilon}}(0)}U^i(x)U^j(x+\frac{y^i-y^j}{\varepsilon})\,dx\\
\leq{}& C\varepsilon^{2+\theta}
\end{align*}
and
\begin{align*}
\int_{\mathbb{R}^2}((V(y^i)-V(a^i))U_{\varepsilon,y^i}^i U_{\varepsilon,y^j}^j\,dx
\leq C\varepsilon^2\int_{\mathbb{R}^2}U^i(x)U^j(x+\frac{y^i-y^j}{\varepsilon})\,dx\leq C\varepsilon^{2+\theta}.
\end{align*}
Therefore we can conclude
\begin{equation}\label{eq4}
I_1=\frac{1}{2}\varepsilon^2\sum_{i=1}^k\int_{\mathbb{R}^2}(V(y^i)-V(a^i))(U^i(x))^2\,dx+O(\varepsilon^{2+\theta}).
\end{equation}
To estimate $I_2$, note that for $a,b>0$ and $p>2$,
\begin{equation}\nonumber
(a+b)^{p}-a^{p}-b^{p}=
O(a^{p-1}b+ab^{p-1})
\end{equation}
and
\begin{equation}\nonumber
(a+b)^{p}-a^{p}-b^{p}-pa^{p-1}b-pab^{p-1}=
\left\{
\begin{aligned}
&O(a^{\frac{p}{2}}b^{\frac{p}{2}}),&\text{if}&\ p\leq3,\\
&O(a^{p-2}b^2+a^{2}b^{p-2}),&\text{if}&\ p>3.
\end{aligned}
\right.
\end{equation}
Then we obtain
\begin{equation}\label{eq5}\nonumber
\begin{aligned}
\int_{\mathbb{R}^{2}}\Big(\sum_{i=1}^k U_{\varepsilon,y^i}^i\Big)^p\,dx&=\sum_{i=1}^k\int_{\mathbb{R}^{2}}(U_{\varepsilon,y^i}^i)^p\,dx+p\sum_{i\neq j}\int_{\mathbb{R}^{2}}(U_{\varepsilon,y^i}^i)^{p-1}U_{\varepsilon,y^j}^j\,dx+O(e^{-\frac{\tau}{\varepsilon}})\\
\end{aligned}
\end{equation}
for some $\tau>0$, which implies that
\begin{equation}\label{eq6}
\begin{aligned}
I_2&=\Big(\frac{1}{2}-\frac{1}{p}\Big)\sum_{i=1}^k\int_{\mathbb{R}^{2}}(U_{\varepsilon,y^i}^i)^p\,dx-\frac{1}{2}\sum_{i\neq j}\int_{\mathbb{R}^2}(U_{\varepsilon,y^i}^i)^{p-1} U_{\varepsilon,y^j}^j\,dx+O(e^{-\frac{\tau}{\varepsilon}})\\
&=\Big(\frac{1}{2}-\frac{1}{p}\Big)\varepsilon^2\sum_{i=1}^k
\int_{\mathbb{R}^2}(U^i(x))^p\,dx+O(e^{-\frac{\tau}{\varepsilon}}).
\end{aligned}
\end{equation}
Hence \eqref{eq2} follows from \eqref{eq7}-\eqref{eq6}.
\end{proof}

\section{Local Pohozaev identities}\label{secb}
In this section, we present the following crucial local Pohozaev identities which are used before
and interesting independently.

\begin{prop}
Suppose that $(u,A_0,A_1,A_2)$ is a weak solution of \eqref{eq1}. Then we have the following local Pohozaev identity:
\begin{equation}\label{P}
\begin{aligned}
\frac{1}{2}\int_{\Omega}\frac{\partial V}{\partial x_k}u^2\,dx={}& \frac{1}{2}\int_{\partial\Omega}(\varepsilon^2|\nabla u|^2+V(x)u^2)\nu_k\,d\sigma-\int_{\partial\Omega}\varepsilon^2\frac{\partial u}{\partial\nu}\frac{\partial u}{\partial x_k}\,d\sigma-\frac{1}{p}\int_{\partial\Omega}|u|^p\nu_k\,d\sigma\\
{}&-\int_{\partial\Omega}(A_1\nu_1+A_2\nu_2)A_ku^2\,d\sigma+\frac{1}{2}\int_{\partial\Omega}
(A_0+A_1^2+A_2^2)u^2\nu_k\,d\sigma,
\end{aligned}
\end{equation}
where $\Omega$ is a bounded open domain of $\mathbb{R}^2$, $k=1,2,$ and $\nu=(\nu_1,\nu_2)$ is the outward unit normal of $\partial\Omega$.
\end{prop}
\noindent
\begin{proof}
To prove them, we will take full advantage of the relationship between $A_0$, $A_1$ and $A_2.$
Assume that $(u,A_0,A_1,A_2)$ is a weak solution of \eqref{eq1}. Multiplying the first equation of \eqref{eq1} by $\overline{\partial_k u+iA_k u}$ and integrating on $\Omega$, we have
\begin{equation}\nonumber
\int_{\Omega}\Big(-\varepsilon^2\Delta u+(V(x)+A_0+A_1^2+A_2^2)u\Big)\overline{\partial_k u+iA_k u}\,dx=
\int_{\Omega}|u|^{p-2}u\overline{\partial_k u+iA_k u}\,dx.
\end{equation}
Observing that
\begin{equation}\nonumber
\Re\Big\{\int_{\Omega}-\varepsilon^2\Delta u\overline{\partial_k u+iA_k u}\,dx\Big\}
=-\varepsilon^2\int_{\Omega}\Delta u\frac{\partial u}{\partial x_k}\,dx
=\frac{1}{2}\int_{\partial\Omega}\varepsilon^2|\nabla u|^2\nu_k\,d\sigma
-\int_{\partial\Omega}\varepsilon^2\frac{\partial u}{\partial\nu}\frac{\partial u}{\partial x_k}\,d\sigma,
\end{equation}
\begin{equation}\nonumber
\begin{aligned}
{}&\Re\Big\{\int_{\Omega}(V(x)+A_0+A_1^2+A_2^2)u\overline{\partial_k u+iA_k u}\,dx\Big\}=\int_{\Omega}(V(x)+A_0+A_1^2+A_2^2)u
\frac{\partial u}{\partial x_k}\,dx\\
={}&-\frac{1}{2}\int_{\Omega}\frac{\partial}{\partial x_k}(V(x)+A_0+A_1^2+A_2^2)u^2\,dx
+\frac{1}{2}\int_{\partial\Omega}(V(x)+A_0+A_1^2+A_2^2)u^2\nu_k\,d\sigma,
\end{aligned}
\end{equation}
and
\begin{equation}\nonumber
\Re\Big\{\int_{\Omega}|u|^{p-2}u\overline{\partial_k u+iA_k u}\,dx\Big\}=\int_{\Omega}|u|^{p-2}u\frac{\partial u}{\partial x_k}\,dx=\frac{1}{p}\int_{\partial\Omega}|u|^p\nu_k\,d\sigma,
\end{equation}
 we have
\begin{equation}\label{aa}
\begin{aligned}
{}&\frac{1}{2}\int_{\partial\Omega}\varepsilon^2|\nabla u|^2\nu_k\,d\sigma
-\int_{\partial\Omega}\varepsilon^2\frac{\partial u}{\partial\nu}\frac{\partial u}{\partial x_k}\,d\sigma
-\frac{1}{2}\int_{\Omega}\frac{\partial}{\partial x_k}(V(x)+A_0+A_1^2+A_2^2)u^2\,dx\\
{}&+\frac{1}{2}\int_{\partial\Omega}(V(x)+A_0+A_1^2+A_2^2)u^2\nu_k\,d\sigma
-\frac{1}{p}\int_{\partial\Omega}|u|^p\nu_k\,d\sigma=0.
\end{aligned}
\end{equation}
Using the Coulomb condition and \eqref{eq1}, we find that
\begin{equation}\nonumber
\begin{aligned}
0&=\partial_2\partial_1A_0-\partial_1\partial_2A_0=\partial_2(A_2u^2)+\partial_1(A_1u^2)\\
&=2u(A_1\partial_1u+A_2\partial_2u)+u^2(\partial_1A_1+\partial_2A_2)\\
&=2u(A_1\partial_1u+A_2\partial_2u).
\end{aligned}
\end{equation}
This shows that $A_1\partial_1u+A_2\partial_2u=0$. Then we can verify that for $k=1,2,$
\begin{equation}\label{bb}
\begin{aligned}
0={}&\Re\Big\{i\int_{\Omega}\overline{\partial_k u+iA_k u}(A_1\partial_1u+A_2\partial_2u)\,dx\Big\}
=\int_{\Omega}(A_kA_1u\partial_1u+A_kA_2u\partial_2u)\,dx\\
={}&-\int_{\Omega}\partial_1(A_kA_1u)u\,dx+\int_{\partial\Omega}A_kA_1u^2\nu_1\,d\sigma
-\int_{\Omega}\partial_2(A_kA_2u)u\,dx+\int_{\partial\Omega}A_kA_2u^2\nu_2\,d\sigma\\
={}&-\int_{\Omega}\partial_1(A_kA_1)u^2\,dx-\int_{\Omega}A_kA_1u\partial_1u\,dx+\int_{\partial\Omega}A_kA_1u^2\nu_1\,d\sigma\\
{}&-\int_{\Omega}\partial_2(A_kA_2)u^2\,dx-\int_{\Omega}A_kA_2u\partial_2u\,dx+\int_{\partial\Omega}A_kA_2u^2\nu_2\,d\sigma\\
={}&-\frac{1}{2}\int_{\Omega}\partial_1(A_kA_1)u^2\,dx
-\frac{1}{2}\int_{\Omega}\partial_2(A_kA_2)u^2\,dx
+\frac{1}{2}\int_{\partial\Omega}A_kA_1u^2\nu_1\,d\sigma
+\frac{1}{2}\int_{\partial\Omega}A_kA_2u^2\nu_2\,d\sigma.
\end{aligned}
\end{equation}
Applying \eqref{bb} with $k=1$ and the equation \eqref{eq1}, we have
\begin{equation}\label{cc}
\begin{aligned}
{}&\int_{\Omega}\frac{\partial}{\partial x_1}(A_0+A_1^2+A_2^2)u^2\,dx
=\int_{\Omega}(\partial_1A_0u^2+\partial_1(A_1^2)u^2+\partial_1(A_2^2)u^2)\,dx\\
={}&\int_{\Omega}(\partial_1A_0u^2-\partial_2(A_1A_2)u^2+\partial_1(A_2^2)u^2)\,dx
+\int_{\partial\Omega}(A_1^2u^2\nu_1+A_1A_2u^2\nu_2)\,d\sigma\\
={}&\int_{\Omega}(\partial_1A_0u^2+(A_2\partial_1A_2-A_1\partial_2A_2)u^2+A_2(\partial_1A_2-\partial_2A_1)u^2)\,dx
+\int_{\partial\Omega}(A_1^2u^2\nu_1+A_1A_2u^2\nu_2)\,d\sigma\\
={}&\int_{\Omega}(\partial_1A_0u^2+(A_2\partial_1A_2+A_1\partial_1A_1)u^2-\frac{1}{2}u^2\partial_1A_0)\,dx
+\int_{\partial\Omega}(A_1^2u^2\nu_1+A_1A_2u^2\nu_2)\,d\sigma\\
={}&\frac{1}{2}\int_{\Omega}(\partial_1A_0u^2+\partial_1(A_1^2)u^2+\partial_1(A_2^2)u^2)\,dx
+\int_{\partial\Omega}(A_1^2u^2\nu_1+A_1A_2u^2\nu_2)\,d\sigma\\
={}&2\int_{\partial\Omega}(A_1^2u^2\nu_1+A_1A_2u^2\nu_2)\,d\sigma.
\end{aligned}
\end{equation}
Similarly, we have
\begin{equation}\label{dd}
\int_{\Omega}\frac{\partial}{\partial x_2}(A_0+A_1^2+A_2^2)u^2\,dx
=2\int_{\partial\Omega}(A_1A_2u^2\nu_1+A_2^2u^2\nu_2)\,d\sigma.
\end{equation}
By substituting \eqref{cc} and \eqref{dd} into \eqref{aa}, we finish the proof.
\end{proof}

\textbf{Acknowledgment.}
This paper was supported by NSFC grants (No.12071169; No. 12226324).

\end{document}